\documentclass[a4paper,11pt]{amsart}

\usepackage{amssymb,amsthm,amsmath,amscd,amsfonts,bbm,mathrsfs}
\usepackage{stmaryrd}
\usepackage{hyperref}
\hypersetup{
  bookmarksnumbered=true,
  colorlinks=true,
  linkcolor=blue,
  citecolor=blue,
  filecolor=blue,
  menucolor=blue,
  urlcolor=blue,
  bookmarksopen=true,
  bookmarksdepth=2,
  pageanchor=true}

\usepackage[cmtip,all]{xy}

\theoremstyle{plain}
\newtheorem{Lemma}[equation]{Lemma}
\newtheorem{Thm}[equation]{Theorem}
\newtheorem{Prop}[equation]{Proposition}
\newtheorem{Cor}[equation]{Corollary}
\theoremstyle{definition}
\newtheorem{Defn}[equation]{Definition}
\newtheorem{Ass}[equation]{Assumption}
\theoremstyle{remark}
\newtheorem{Remark}[equation]{Remark}
\newtheorem*{Remark*}{Remark}

\newcommand{\bLe}{\begin{Lemma}}
\newcommand{\eLe}{\end{Lemma}}
\newcommand{\bTh}{\begin{Thm}}
\newcommand{\eTh}{\end{Thm}}
\newcommand{\bPr}{\begin{Prop}}
\newcommand{\ePr}{\end{Prop}}
\newcommand{\bCo}{\begin{Cor}}
\newcommand{\eCo}{\end{Cor}}
\newcommand{\bDe}{\begin{Defn}}
\newcommand{\eDe}{\end{Defn}}
\newcommand{\bAs}{\begin{Ass}}
\newcommand{\eAs}{\end{Ass}}
\newcommand{\bReX}{\begin{Remark*}}
\newcommand{\eReX}{\end{Remark*}}
\newcommand{\bRe}{\begin{Remark}}
\newcommand{\eRe}{\end{Remark}}
\newcommand{\bproof}{\begin{proof}}
\newcommand{\eproof}{\end{proof}}

\numberwithin{equation}{section}

\newcommand{\EEE}{\mathcal{E}}

\newcommand{\KKK}{\mathcal{K}}
\newcommand{\LLL}{\mathcal{L}}
\newcommand{\MMM}{\mathcal{M}}
\newcommand{\NNN}{\mathcal{N}}
\newcommand{\OOO}{\mathcal{O}}
\newcommand{\PPP}{\mathcal{P}}
\newcommand{\QQQ}{\mathcal{Q}}

\newcommand{\SSS}{\mathcal{S}}
\newcommand{\TTT}{\mathcal{T}}

\newcommand{\XXX}{\mathcal{X}}

\newcommand{\Fp}{\mathfrak{p}}
\newcommand{\Fq}{\mathfrak{q}}

\newcommand{\CC}{{\mathbb{C}}}

\newcommand{\RR}{{\mathbb{R}}}

\newcommand{\ZZ}{{\mathbb{Z}}}

\newcommand{\One}{{\mathbbm{1}}}

\newcommand{\mmod}{{\operatorname{-mod}}}
\newcommand{\MMod}{{\operatorname{-Mod}}}
\newcommand{\pperf}{{\operatorname{-perf}}}
\newcommand{\pproj}{{\operatorname{-proj}}}
\newcommand{\PProj}{{\operatorname{-Proj}}}

\DeclareMathOperator{\cone}{cone}
\DeclareMathOperator{\even}{even}
\DeclareMathOperator{\et}{et}
\DeclareMathOperator{\id}{id}
\DeclareMathOperator{\ind}{ind}
\DeclareMathOperator{\liset}{lis-et}
\DeclareMathOperator{\lat}{lat}

\DeclareMathOperator{\qc}{qc}

\DeclareMathOperator{\red}{red}
\DeclareMathOperator{\res}{res}

\DeclareMathOperator{\stmod}{stmod}
\DeclareMathOperator{\supp}{supp}

\DeclareMathOperator{\Aut}{Aut}

\DeclareMathOperator{\Cohuh}{Coh-uh}
\DeclareMathOperator{\End}{End}

\DeclareMathOperator{\Frac}{Frac}

\DeclareMathOperator{\Hom}{Hom}

\DeclareMathOperator{\Ker}{Ker}
\DeclareMathOperator{\Max}{Max}
\DeclareMathOperator{\Proj}{Proj}
\DeclareMathOperator{\Perf}{Perf}

\DeclareMathOperator{\Spc}{Spc}
\DeclareMathOperator{\Spcsurj}{Spc-surj}

\DeclareMathOperator{\StMod}{StMod}
\DeclareMathOperator{\Spec}{Spec}

\begin{document}

\title{The Balmer spectrum of certain Deligne-Mumford stacks}
\author{Eike Lau}
\date{\today}
\address{Eike Lau, Fakult\"{a}t f\"{u}r Mathematik,
Universit\"{a}t Bielefeld, D-33501 Bielefeld}

\begin{abstract}
We consider a Deligne--Mumford stack $X$ which is the quotient of an affine scheme $\operatorname{Spec}A$ by the action of a finite group $G$ and show that the Balmer spectrum of the tensor triangulated category of perfect complexes on $X$ is homeomorphic to the space of homogeneous prime ideals in the group cohomology ring $H^*(G,A)$. 
\end{abstract}

\maketitle

\setcounter{tocdepth}{1}
\tableofcontents

\section{Introduction}

Let $\XXX$ be an algebraic stack. The category $\Perf(\XXX)$ of perfect complexes of $\OOO_\XXX$-modules is a tensor triangulated category. One can ask for a classification of all thick tensor ideals in this category, or equivalently for a description of the Balmer spectrum
\[
\Spc(\Perf(\XXX)),
\]
the space of all prime thick tensor ideals in $\Perf(\XXX)$, as defined in \cite{Balmer:Spectrum}. An answer is known at least in the following cases. 

If $\XXX$ is a qcqs scheme, the thick tensor ideals in $\Perf(\XXX)$ are classified in terms of their support by Thomason \cite{Thomason:Classification}. This yields a homeomorphism
\begin{equation}
\label{Eq:PerfX-X}
\Spc(\Perf(\XXX))\cong|\XXX|,
\end{equation}
where $|\XXX|$ is the underlying topological space of $\XXX$, by Balmer 
\cite[Thm.\ 4.1]{Balmer:Guide}, which relies on \cite{Buan-Krause-Solberg} when $X$ is not noetherian. More generally, a homeomorphism \eqref{Eq:PerfX-X} exists when $\XXX$ is a quasicompact algebraic stack with separated diagonal and $\XXX$ is tame, which means that the geometric stabilizer groups of $\XXX$ are finite and geometrically reductive, by Hall \cite{Hall:Tame}.

In a different direction, if $\XXX=BG$ is the classifying space of a finite group $G$ over a field $k$, then $\Perf(\XXX)$ is equivalent to $D^b(kG\mmod)$, the bounded derived category of finite modules over the group algebra $kG$. In this case, thick tensor ideals in the stable module category of $kG$ are classified in terms of their cohomological support by Benson--Carlson--Rickard \cite{BCR:Thick} and Benson--Iyengar--Krause \cite{BIK:Stratifying}; see also \cite{Carlson-Iyengar:Thick} for a direct proof of a non-stable version. These results yield a homeomorphism
\begin{equation}
\label{Eq:PerfX-RkG}
\Spc(\Perf(\XXX))\cong\Spec^h(R_{k,G})
\end{equation}
where $R_{k,G}=H^*(G,k)$ is the cohomology ring and $\Spec^h$ is the space of all homogeneous prime ideals, by Balmer \cite[Prop.\ 8.5]{Balmer:Spectra3}. The homeomorphism \eqref{Eq:PerfX-RkG} is complementary to \eqref{Eq:PerfX-X} because the stack $\XXX=BG$ is tame iff the characteristic of $k$ does not divide the order of $G$, in which case $\Spc(\Perf(\XXX))$ is just a point. More generally, if $G$ is a finite group scheme over $k$, a similar classification of the thick tensor ideals in the stable module category of $kG$ and hence the homeomorphism \eqref{Eq:PerfX-RkG} for $\XXX=BG$ exist by Benson--Iyengar--Krause--Pevtsova \cite{BIKP:Stratification}.

\subsection*{Our main example}
We consider a quotient stack $\XXX=[\Spec(A)/G]$ where $G$ is a finite group that acts on a commutative ring $A$. Let $\TTT_{A,G}=\Perf(\XXX)$. This category is equivalent to the category $D^b(AG)_{A\pperf}$ of bounded complexes of $AG$-modules which are perfect as complexes of $A$-modules, viewed as a full subcategory of the derived category $D(AG)$. By \cite{Balmer:Spectra3} there is a natural continuous comparison map
\begin{equation*}
\rho_{A,G}:\Spc(\TTT_{A,G})\to\Spec^h(R_{A,G})
\end{equation*}
with $R_{A,G}=H^*(G,A)$. If $A$ is a field with the trivial action of $G$, which we will call the punctual case, then $\rho_{A,G}$ is the homeomorphism \eqref{Eq:PerfX-RkG}. If $G$ is the trivial group, then $\rho_{A,G}$ is the homeomorphism \eqref{Eq:PerfX-X} for the affine scheme $\XXX=\Spec A$. The following is our main result.

\bTh
\label{Th:main-intro}
The map $\rho_{A,G}$ is a homeomorphism in all cases.
\eTh

Even the case $A=\ZZ$ of integral representations seems to be new.

\subsection*{A stable version}
The category $\Perf(AG)$ of perfect complexes of $AG$-modules is a tensor ideal in $\TTT_{A,G}$, and hence the Verdier quotient 
\[
\SSS_{A,G}=\TTT_{A,G}/\Perf(AG)
\] 
is a tensor triangulated category again. As observed in \cite{Barthel:StratifyingIntegral}, the following variant of Theorem \ref{Th:main-intro} is an immediate consequence.

\bCo
\label{Co:main-intro}
The map $\rho_{A,G}$ induces a homeomorphism
\[
\Spc(\SSS_{A,G})\cong\Proj(R_{A,G}).
\]
\eCo

The category $\SSS_{A,G}$ can be viewed as the stable category of the Frobenius category $\lat(A,G)$ of all finite $A$-projective $AG$-modules, so Corollary \ref{Co:main-intro} is in line with the classical results of modular representation theory mentioned above. If $G$ acts trivially on $A$, another stable category of $AG$-modules with different behaviour was introduced in \cite{BIK:ModuleCategories}, using the Frobenius category of all finitely presented $AG$-modules with $A$-split exact sequences.

\subsection*{Recent literature}

In the first version of this article, Theorem \ref{Th:main-intro} was proved only when the ring $A$ is regular. Afterwords, related results on stratification of compactly generated triangulated categories of $AG$-modules in the case where $A$ is regular with trivial $G$-action appeared in \cite{Barthel:StratifyingIntegral,Barthel:StratifyingViaHomotopy} and in \cite{BIKP:Fibrewise}. The fact that Theorem \ref{Th:main-intro} implies Corollary \ref{Co:main-intro} was observed in \cite{Barthel:StratifyingIntegral}.

\subsection*{Strategy}
The proof of Theorem \ref{Th:main-intro} is based on a reduction to the punctual case and the case of affine schemes along the following lines.  If $A$ is a field with possibly non-trivial action of $G$, the result follows from the punctual case by a form of Galois descent. In general we can assume that $A$ is essentially of finite type over $\ZZ$; then it will be sufficient to show that $\rho_{A,G}$ is bijective.
There is a commutative diagram of continuous maps
\begin{equation}
\label{Dia:fiber-intro}
\vcenter{
\xymatrix@M+0.2em@C+1em{
\Spc(\TTT_{A,G}) \ar[r]^-{\rho_{A,G}} \ar[d]^{\pi_\TTT} &  \Spec^h(R_{A,G}) \ar[d]^{\pi_R} \\
\Spc(\TTT_{A^G}) \ar[r]^\sim_{\rho_{A^G}} & \Spec (A^G)
}}
\end{equation}
where $A^G$ is the ring of $G$-invariants and $\TTT_{A^G}=\Perf(\Spec(A^G))$ is the associated category of perfect complexes; note that $A^G$ is the degree zero component of the graded ring $R_{A,G}$. Here $\rho_{A^G}$ is a homeomorphism by the case of affine schemes. Hence the map $\rho_{A,G}$ is bijective iff it restricts to a bijective map between each fiber of $\pi_\TTT$ and the corresponding fiber of $\pi_R$. This map between the fibers will be related with the field case of Theorem \ref{Th:main-intro} as follows. For a prime ideal $\Fq\in\Spec(A^G)$ we consider the reduced fiber of the ring $A$ over $\Fq$,
\[
A(\Fq)=(A\otimes_{A^G}k(\Fq))_{\red}.
\]
Explicitly this means that $A(\Fq)=k(\Fp_1)\times\ldots\times k(\Fp_r)$ where $\Fp_1,\ldots,\Fp_r$ are the prime ideals of $A$ over $\Fq$. Functoriality gives the following commutative diagram, where the index $\Fq$ means fiber over $\Fq$ in \eqref{Dia:fiber-intro}.
\begin{equation*}
\vcenter{
\xymatrix@M+0.2em@C+1em{
\Spc(\TTT_{A(\Fq),G}) \ar[r]^-{\rho_{A(\Fq),G}} 
\ar[d]_{j_\TTT} 
& \Spec^h(R_{A(\Fq),G}) 
\ar[d]^{j_R} 
\\
\Spc(\TTT_{A,G})_\Fq \ar[r]^-{(\rho_{A,G})_\Fq} & \Spec^h(R_{A,G})_\Fq.
}}
\end{equation*}
The field case of Theorem \ref{Th:main-intro} implies that $\rho_{A(\Fq),G}$ is a homeomorphism. We will show that $j_R$ is a homeomorphism and that $j_\TTT$ is surjective, at least when $A$ is noetherian. It follows that $(\rho_{A,G})_\Fq$ is bijective as required.

\subsection*{Change of coefficients}

The preceding assertions on the maps $j_R$ and $j_\TTT$ fall under the question how the spaces $\Spec^h(R_{A,G})$ and $\Spc(\TTT_{A,G})$ vary with the ring $A$. A $G$-equivariant homomorphism of commutative rings $A\to B$ gives rise to a base change map for the spectra of the cohomology rings
\begin{equation}
\label{Eq:bc-Spech-intro}
\Spec^h(R_{B,G})\to
\Spec^h(R_{A,G})\times_{\Spec(A^G)}\Spec(B^G)
\end{equation}
and a base change map for the Balmer spectra
\begin{equation}
\label{Eq:bc-Spc-intro}
\Spc(\TTT_{B,G})\to\Spc(\TTT_{A,G})\times_{\Spec(A^G)}\Spec(B^G).
\end{equation}
In the case $B=A(\Fq)$ these maps can be identified with $j_R$ and $j_\TTT$. The following general result on base change in group cohomology shows in particular that $j_R$ is a homeomorphism.

\bTh
\label{Th:univ-homeo-intro}
The base change map \eqref{Eq:bc-Spech-intro} is a homeomorphism
if $B$ is a localisation of a quotient of $A$.
\eTh

See Corollary \ref{Co:Spech-Spec-cart} and the underlying Theorem \ref{Th:univ-homeo}, which says that the natural ring homomorphism $(R_{A,G})_{\even}\otimes_{A^G}B^G\to(R_{B,G})_{\even}$ induces a universal homeomorphism of affine schemes if $B$ is a localisation of a quotient of $A$. The proof of that result proceeds by a reduction to three basic cases: the case where $B$ is a localisation of $A$, the case $B=A/I$ for a nilpotent ideal $I$, and the case $B=A/tA$ for an $A$-regular element $t\in A^G$, under a noetherian assumption. 

With hindsight, Theorem \ref{Th:main-intro} identifies \eqref{Eq:bc-Spech-intro} and \eqref{Eq:bc-Spc-intro}. A direct proof of the analogue of Theorem \ref{Th:univ-homeo-intro} for the map \eqref{Eq:bc-Spc-intro} seems to be difficult, but in the above three basic cases we can at least show that \eqref{Eq:bc-Spc-intro} is surjective, using a general surjectivity criterion of Balmer \cite{Balmer:Surjectivity}, saying that a functor between tensor triangulated categories which detects tensor nilpotence induces a surjective map on the Balmer spectra. By formal arguments one deduces that \eqref{Eq:bc-Spc-intro} is surjective when $A$ is noetherian and $B=A(\Fq)$, in other words the map $j_\TTT$ is surjective in the noetherian case.

\subsection*{Structure of the article}

In section \ref{Se:Prime-spectra} we review the comparison map $\rho$ and its relation with cohomological support for general tensor triangulated categories. 
In section \ref{Se:modules-skew} we study the category $D^b(AG)_{A\pperf}$ from a purely algebraic point of view.
In section \ref{Se:perfect} we introduce the category $\TTT_{A,G}$ of perfect complexes on $\XXX$ and show that it is equivalent to the algebraic category of section \ref{Se:modules-skew}. 
In section \ref{Se:comparison-basic} we study basic properties of the comparison map for these categories.
The field case of Theorem \ref{Th:main-intro} is established in section \ref{Se:field}. 
Theorem \ref{Th:univ-homeo-intro} is proved in section \ref{Se:change-coefficient}. 
The corresponding surjectivity results for the map $j_\TTT$ are proved in sections \ref{Se:TensorNil} and \ref{Se:ChangeCoeffSpc}.
Finally, the main result is proved in section \ref{Se:main}, and Corollary \ref{Co:main-intro} is deduced in section \ref{Se:stable}.

\subsection*{Notation} 

For a not necessarily commutative ring $R$ we denote by $R\MMod$ the category of left $R$-modules and by $R\mmod$ the category of finitely generated left $R$-modules. The latter is abelian if $R$ is left noetherian.

\medskip\smallskip
This note was inspired by a talk given by H.\ Krause on \cite{BIKP:Stratification} in Bielefeld. I am grateful to Henning Krause and Bill Crawley-Boevey for helpful comments. I thank the referee for very careful reading of the manuscript and many helpful suggestions for improvement.

\section{Prime spectra and cohomological support}
\label{Se:Prime-spectra}

In this section we recall some aspects of Balmer's theory of prime spectra of tensor triangulated categories. Let $\TTT$ be an essentially small tensor triangulated category with unit object $\One$ and let 
\[
R=R_\TTT=\End^*_\TTT(\One)=\bigoplus_{n\in\ZZ}\Hom_\TTT(\One,\One[n])
\]
as a graded ring. This is a graded-commutative ring by the obvious composition of morphisms, or equivalently by the tensor product of morphisms; see \cite[Prop.~3.3]{Balmer:Spectra3}. By \cite[Thm.~5.3]{Balmer:Spectra3} there is a continuous map, denoted there by $\rho^\bullet$ and called the comparison map in \cite{Balmer:Guide}, 
\begin{equation}
\label{Eq:rho}
\rho=\rho_\TTT:\Spc(\TTT)\to\Spec^h(R)
\end{equation}
from the space of prime thick tensor ideals in $\TTT$ to the space of homogeneous prime ideals in $R$, which is defined by the following property. If $\rho(\PPP)=\Fp$, a homogeneous element $a\in R$ satisfies
$a\not\in\Fp$ iff $\cone(a)\in\PPP$.

If the ring $R$ is noetherian, then $\rho$ is surjective by \cite[Thm.\ 7.3]{Balmer:Spectra3}.

\subsection*{Support}

For a graded $R$-module $M$ we have the support 
\[
\supp_R(M)=\{\Fp\in\Spec^h(R)\mid M_\Fp\ne 0\}
\]
where $\Spec^h(R)$ is the set of homogeneous prime ideals of $R$ and $M_\Fp$ is the homogeneous localisation of $M$ at $\Fp$. For an object $X$ of $\TTT$ there are two notions of support: The canonical support 
\[
\supp(X)=\{\PPP\in\Spc(\TTT)\mid X\not\in\PPP\},
\]
which is a closed subset of $\Spc(\TTT)$, and the cohomological support
\[
V(X)=\supp_R(\End^*_\TTT(X))
\]
as a subset of $\Spec^h(R)$, where $M_X=\End^*_\TTT(X)$ becomes a graded $R$-module by the graded ring homomorphism $R\to M_X$ defined by $f\mapsto f\otimes\id_X$, whose image is graded-central; the same $R$-module structure on $M_X$ arises from the tensor product of morphisms; see \cite[Prop.~3.3]{Balmer:Spectra3}.

The cohomological support $V(X)$ can also be described using localisations of $\TTT$ as in \cite[Constr.~3.5]{Balmer:Spectra3}. For $\Fp\in\Spec^h(R)$, the localisation $\TTT_\Fp$ has the same objects as $\TTT$ and homomorphism groups $\Hom_{\TTT_\Fp}(X,Y)=(\Hom^*_\TTT(X,Y)_\Fp)^0$, where $^0$ means degree zero part. Let $X_\Fp\in\TTT_\Fp$ be the image of $X\in\TTT$. Then
\begin{equation}
\label{Eq:FpinV(X)}
\Fp\in V(X)
\;\Longleftrightarrow\;
(M_X)_\Fp\ne 0
\;\Longleftrightarrow\;
X_\Fp\ne 0.
\end{equation}
The category $\TTT_\Fp$ is equivalent to the Verdier localisation of $\TTT$ with respect to the thick tensor ideal generated by $\cone(a)$ for all homogeneous elements $a\in R\setminus\Fp$ by \cite[Thm.\ 3.6]{Balmer:Spectra3}, in particular $\TTT_\Fp$ is a tensor triangulated category such that the localisation $\TTT\to\TTT_\Fp$ is an exact tensor functor.

\bLe
\label{Le:V-cone-a}
For a homogeneous element $a\in R$ we have
\[
V(\cone(a))=V(a)=\{\Fp\in\Spec^h(R)\mid a\in\Fp\}.
\]
\eLe

\bproof
Using \eqref{Eq:FpinV(X)} and the exactness of $\TTT\to\TTT_\Fp$ we have
$\Fp\in V(\cone(a))$ iff $\cone(a)_\Fp\ne 0$ iff $a$ not invertible in $\TTT_{\Fp}$ iff $a\in\Fp$. 
\eproof

One can ask if $(\Spec^h(R),V)$ is a support data on $\TTT$ as in \cite[Def.\ 3.1]{Balmer:Spectrum}.

\bLe
\label{Le:support-data}
The cohomological support $(\Spec^h(R),V)$ is a support data on $\TTT$ iff for all $X,Y\in\TTT$ the set $V(X)$ is closed and 
\begin{equation}
\label{Eq:Tensor-Product-Formula}
V(X\otimes Y)=V(X)\cap V(Y).
\end{equation}
\eLe

\bproof
The conditions are part of the axioms of support data. Since the localisation $\TTT_p$ is a tensor triangulated category, the remaining axioms of support data are easily verified using \eqref{Eq:FpinV(X)}.
\eproof

\subsection*{The comparison map and support}

\bDe
\label{De:End-finite}
The tensor triangulated category $\TTT$ will be called End-finite if for each $X\in\TTT$ the $R$-module $M_X=\End^*_\TTT(X)$ is noetherian.
\eDe

Clearly $\TTT$ is End-finite iff the ring $R=\End^*_\TTT(\One)$ is noetherian and the $R$-module $M_X=\End^*_\TTT(X)$ is finite for each $X\in\TTT$, moreover this implies that the $R$-module $\Hom^*_\TTT(X,Y)$ is finite for all $X,Y\in\TTT$ since the latter is a direct summand of $\End^*_\TTT(X\oplus Y)$. 

If $\TTT$ is End-finite, $V(X)$ is closed in $\Spec^h(R)$ for each $X\in\TTT$.

\bPr
\label{Pr:supp-V}
For each object $X$ of $\TTT$ we have
\[
\rho(\supp(X))\subseteq V(X),
\]
with equality if\/ $(\Spec^h(R),V)$ is a support data or\/ $\TTT$ is End-finite and rigid.
\ePr

\bproof
Let $\PPP\in\Spc(\TTT)$ and $\Fp=\rho(\PPP)$. We consider the multiplicative sets $\SSS$ in the category $\TTT$ and $S$ in the ring $R$ defined by
\[
\SSS=\{f:X\to Y \text{ in } \TTT\mid\cone(f)\in\PPP\},
\]
\[
S=\SSS\cap R.
\]
Then $S$ is the set of homogeneous elements in $R\setminus\Fp$ by the definition of the map $\rho$. The corresponding localisations of $\TTT$ are related by functors
\[
\TTT\longrightarrow S^{-1}\TTT\longrightarrow\SSS^{-1}\TTT=\TTT/\PPP,
\]
where $S^{-1}\TTT=\TTT_\Fp$, and $\SSS^{-1}\TTT$ is the Verdier localisation at $\PPP$. We get
\[
(M_X)_\Fp=0
\;\Longleftrightarrow\;
X=0 \text{ in } \TTT_\Fp
\;\Longrightarrow\;
X=0 \text{ in } \TTT/\PPP
\;\Longleftrightarrow\;
X\in\PPP,
\]
hence $\PPP\in\supp(X)\Longrightarrow\Fp\in\supp_R(M_X)$. This proves the first assertion.

To prove the second assertion, for a given $\Fp\in\supp_R(M_X)$ we have to find a prime ideal $\PPP\in\supp(X)$ with $\rho(\PPP)=\Fp$. We can replace $\TTT$ by the localisation $\TTT_\Fp$, cf.\ \cite[Thm.\ 5.4]{Balmer:Spectra3}. Then $R$ is a local graded ring with unique graded maximal ideal $\Fp$. Let $\MMM$ be the multiplicative set of all objects of $\TTT$ of the form
\[
X^{\otimes n}\otimes\cone(a_1)\otimes\ldots\otimes\cone(a_r)
\]
with $n\ge 0$ and homogeneous elements $a_1,\ldots,a_r\in\Fp$. We claim that $0\not\in\MMM$. Then by \cite[Lemma 2.2]{Balmer:Spectrum} there is a $\PPP\in\Spc(\TTT)$ with $\MMM\cap\PPP=\emptyset$, in particular $X\not\in \PPP$ and $\cone(a)\not\in\PPP$ for each homogeneous element $a\in\Fp$. The first condition means that $\PPP\in\supp(X)$, the second condition implies that $\Fp\subseteq\rho(\PPP)$ and thus $\Fp=\rho(\PPP)$ since $\Fp$ was assumed to be maximal. So it remains to verify that $0\not\in\MMM$.

If $(\Spec^h(R),V)$ is a support data, the tensor product formula \eqref{Eq:Tensor-Product-Formula} implies that $\Fp\in V(Z)$ for every $Z\in\MMM$ because $\Fp\in V(X)$ and $\Fp\in V(\cone (a))$ for every $a\in\Fp$ by Lemma \ref{Le:V-cone-a}. Hence $0\not\in\MMM$. Assume that $\TTT$ is End-finite and rigid. We have $X\ne 0$ since $\Fp\in V(X)$, hence $X^{\otimes n}\ne 0$ for $n\ge 0$ since $\TTT$ is rigid. So it suffices to verify that $Y\ne 0$ in $\TTT$ implies $Y\otimes\cone(a)\ne 0$ for every homogeneous element $a\in\Fp$. The triangle $Y\xrightarrow a Y\to Y\otimes\cone(a)\to^+$ gives an exact sequence
\[
\End^*_\TTT(Y)\overset a \longrightarrow \End^*_\TTT(Y)\longrightarrow\Hom^*_\TTT(Y,Y\otimes\cone(a)).
\]
Since $\End^*_\TTT(Y)$ is a non-zero finite $R$-module, here the cokernel of $a$ is non-zero by the graded version of Nakayama's Lemma, thus $Y\otimes\cone(a)\ne 0$. 
\eproof

\bCo
\label{Co:rho-bij-homeo}
If\/ $(\Spec^h(R),V)$ is a support data or\/ $\TTT$ is End-finite and rigid,
then $\rho$ is bijective iff $\rho$ is a homeomorphism.
\eCo

\bproof
We have to show: if $\rho$ is bijective and $Z\subseteq \Spc(\TTT)$ is closed, then $\rho(Z)$ is closed. The definition of the topology on $\Spc(\TTT)$ implies that $Z$ is an intersection of sets of the form $\supp(X_i)$ with $X_i\in\TTT$. Proposition \ref{Pr:supp-V} yields $\rho(\supp(X_i))=V(X_i)$, which is closed in $\Spec^h(R)$ by the assumption. Since $\rho$ is injective, $\rho(Z)$ is the intersection of the sets $V(X_i)$ and thus closed.
\eproof

\bRe
One can ask if the conclusion of Corollary \ref{Co:rho-bij-homeo} holds for general tensor triangular categories; see \cite[Lemma 2.1]{DellAmbrogio-Stanley:Affine} and \cite{DellAmbrogio-Stanley:Affine-Err} for a discussion.
\eRe

\bPr
\label{Pr:rho-homeo-supp}
If the category $\TTT$ is rigid and $\rho$ is a homeomorphism, then 
$
\rho(\supp(X))=V(X)
$
for each $X\in\TTT$.
\ePr

\bproof
Let $\PPP\in\Spc(\TTT)$ and $\Fp=\rho(\PPP)$. As in the proof of Proposition \ref{Pr:supp-V} we consider the natural functor $j:\TTT_\Fp\to\TTT/\PPP$ from the localisation at $\Fp$ to the Verdier quotient by $\PPP$. We have to show that $X=0$ in $\TTT_\Fp$ iff $X=0$ in $\TTT/\PPP$. In fact we show that $j$ is an equivalence. By \cite[Thm.~5.4]{Balmer:Spectra3} there is a cartesian diagram of topological spaces
\begin{equation*}
\vcenter{
\xymatrix@M+0.2em@C+1em{
\Spc(\TTT_\Fp) \ar[r]^-{\rho_{\TTT_{\Fp}}} \ar[d]_{\Spc(q)} &
\Spec^h(R_\Fp) \ar[d] \\
\Spc(\TTT) \ar[r]^-{\rho_{\TTT}} &
\Spec^h(R) 
}}
\end{equation*}
where $q:\TTT\to\TTT_\Fp$ is the localisation functor. Since $\rho_\TTT$ is a homeomorphism, $\rho_{\TTT_\Fp}$ is a homeomorphism. The space $\Spec^h(R_\Fp)$ has a unique closed point which maps to $\Fp$ in $\Spec^h(R)$. Hence $\Spc(\TTT_\Fp)$ has a unique closed point $\PPP'$ which maps to $\PPP$ in $\Spc(\TTT)$.  Since $\TTT_\Fp$ is rigid, $\PPP'$ is the zero ideal of $\TTT_\Fp$; see \cite[Prop.~4.2]{Balmer:Spectra3}. Explicitly this shows that $\PPP=q^{-1}(0)$, in particular $q$ maps $\PPP$ to zero. Hence $q$ induces a functor $\TTT/\PPP\to\TTT_\Fp$ which is an inverse of $j$.
\eproof

\subsection*{Products and colimits}

Let us record how the Balmer spectrum and the comparison map behave under finite products and filtered 2-colimits of $\TTT$.

\bLe
\label{Le:Spc-prod}
Let $\TTT=\TTT_1\times\TTT_2$ with tensor triangulated categories $\TTT_i$ and $R_i=\End^*_{\TTT_i}(\One)$. Then there are decompositions into open subspaces 
\[
\Spc(\TTT)=\Spc(\TTT_1)\sqcup\Spc(\TTT_1)
\quad\text{and}\quad
\Spec^h(R)=\Spec^h(R_1)\sqcup\Spec^h(R_2)
\] 
such that $\rho_\TTT=\rho_{\TTT_1}\sqcup\rho_{\TTT_2}$.
\eLe

\bproof
This is straightforward; note that $R=R_1\times R_2$. See \cite[Ch.~1, Ex.~22]{AM:Commutative} for the corresponding assertion for the spectrum of non-graded rings.
\eproof

\bLe
\label{Le:Spc-colim}
Let $\TTT=\varinjlim_i\TTT_i$ be a filtered $2$-colimit of tensor triangulated categories $\TTT_i$ and $R_i=\End^*_{\TTT_i}(\One)$. Then we have
\[
\Spc(\TTT)=\varprojlim_i\Spc(\TTT_i)
\qquad\text{and}\qquad
\Spec^h(R)=\varprojlim_i\Spec^h(R_i)
\]
as topological spaces such that $\rho_\TTT=\varprojlim_i\rho_{\TTT_i}$.
\eLe

\bproof
This is straightforward. See \cite[\href{https://stacks.math.columbia.edu/tag/078L}{Exercise 078L}]{stacks-project} and \cite[Prop.~10.53]{GW:AG1} for the corresponding assertion for the spectrum of non-graded rings.
\eproof

\section{Modules over the skew group ring}
\label{Se:modules-skew}

Let $A$ be a commutative ring with an action of a finite group $G$. We denote by $AG$ the skew group ring for this action, so $AG$ is the free $A$-module with basis $G$, and the multiplication in $AG$ is defined by $(ag)(bh)=ag(b)gh$ for $a,b\in A$ and $g,h\in G$. Let $D(AG)=D(AG\MMod)$ be the derived category of left $AG$-modules. The bounded above derived category $D^-(AG)$ is equivalent to the homotopy category $K^-(AG\PProj)$ of bounded above complexes of projective $AG$-modules, which carries the tensor product
\begin{equation}
\label{Eq:PxQPxAQ}
P\otimes Q=P\otimes_{A}Q
\end{equation}
with the diagonal action of $G$. This makes $D^-(AG)$ into a tensor triangulated category. We denote by 
\[
D^b(AG)_{A\pproj}
\] 
the full subcategory of $D(AG)$ whose objects are the bounded complexes of $AG$-modules which are finite projective $A$-modules, and by
\[
D^b(AG)_{A\pperf}
\]
the full subcategory of $D(AG)$ whose objects are the complexes of $AG$-modules which are $A$-perfect, i.e.\ isomorphic in $D(A)$ to a bounded complex of finite projective $A$-modules. We will verify that these two categories are equivalent. It is easy to see that they are triangulated subcategories of $D^-(AG)$ and that (at least) $D^b(AG)_{A\pperf}$ is stable under the tensor product. 

\bLe
\label{Le:DAGAperf-fin-res}
For every $P\in D^b(AG)_{A\pperf}$ there is a bounded above complex of finite projective $AG$-modules with a quasi-isomorphism $P'\to P$.
\eLe

\bproof
This is standard; see for example \cite[\href{https://stacks.math.columbia.edu/tag/064Z}{Lemma 064Z}]{stacks-project}, where all rings are assumed to be commutative, but that assumption is not used in the proof. Let us sketch a direct argument. Assume $\alpha:Q\to P$ is chosen where $Q=[Q^m\to\ldots\to Q^n]$ is a finite complex of finite projective $AG$-modules such that $H^i(\alpha)$ is surjective for $i=m$ and bijective for $i>m$; equivalently $H^i(\cone(\alpha))=0$ for $i\ge m$. Since $\cone(\alpha)$ is $A$-perfect, it follows that $H^{m-1}(\cone(\alpha))$ is a finite $A$-module, so we can choose a surjective homomorphism of $AG$-modules $Q^{m-1}\to H^{m-1}(\cone(\alpha))$ where $Q^{m-1}$ is finite projective. This gives an extension $Q'=[Q^{m-1}\to Q^m\to\ldots\to Q^n]$ of $Q$ and an extension $\alpha':Q'\to P$ of $\alpha$ such that $H^i(\cone(\alpha'))=0$ for $i\ge m-1$. By infinite repetition the growing $Q$ gives $P'$.
\eproof

\bLe
\label{Le:DAG-Aproj-Aperf}
The inclusion $D^b(AG)_{A\pproj}\to D^b(AG)_{A\pperf}$ is an equivalence of triangulated categories. The resulting tensor product on $D^b(AG)_{A\pproj}$ is given by \eqref{Eq:PxQPxAQ}.
\eLe

\bproof
For $P\in D^b(AG)_{A\pperf}$ let $P'\to P$ be a quasi-isomorphism as in Lemma \ref{Le:DAGAperf-fin-res}. For sufficiently small $n\in\ZZ$, the truncation $P''=\tau_{\ge n}P'$ is quasi-isomorphic to $P$ as well. Then $P''$ lies in $D^b(AG)_{A\pproj}$ because $P''$ is $A$-perfect and bounded with finite $A$-projective components except possibly in the left-most degree $n$. This proves the first assertion. For the second assertion we note that for quasi-isomorphisms $P'\to P$ and $Q'\to Q$ with $P',Q'\in K^-(AG\PProj)$ and $P,Q\in D^b(AG)_{A\pproj}$ the resulting homomorphism $P'\otimes_A Q'\to P\otimes_A Q$ is a quasi-isomorphism again.
\eproof

\bLe
There is an isomorphism of graded rings
\begin{equation}
\label{Eq:End*DAGA-H*GA}
\End^*_{D(AG)}(A)\cong H^*(G,A),
\end{equation}
and for two $AG$-modules $M$ and $N$ such that $M$ is $A$-projective there is a natural isomorphism of graded modules with respect to \eqref{Eq:End*DAGA-H*GA},
\begin{equation}
\label{Eq:Hom*DAGMN-H*GM}
\Hom^*_{D(AG)}(M,N)\cong H^*(G,\Hom_A(M,N)).
\end{equation}
\eLe

Again this is standard, at least when $G$ acts trivially on $A$, but this condition is not essential. We sketch a proof for completeness.

\bproof
We have $\Hom^*_{D(AG)}(M,N)=\Hom^*_{D(AG)}(A,\Hom_A(M,N))$, so in \eqref{Eq:Hom*DAGMN-H*GM} it suffices to treat the case $M=A$. If $P'\to\ZZ$ is a $\ZZ G$-projective resolution, then $P=A\otimes_\ZZ P'\to A$ is an $AG$-projective resolution, and the complex $\Hom_{AG}(P,N)$ is isomorphic to $\Hom_{\ZZ G}(P',N)$. This gives \eqref{Eq:Hom*DAGMN-H*GM} and \eqref{Eq:End*DAGA-H*GA} as graded abelian groups. The cup product on $H^*(G,A)$ corresponds to the tensor product in $\End^*_{D(AG)}(A)$, which coincides with the composition product by \cite[Prop.~3.3]{Balmer:Spectra3}, and hence \eqref{Eq:End*DAGA-H*GA} is an isomorphism of graded rings. Similarly \eqref{Eq:Hom*DAGMN-H*GM} is an isomorphism of graded modules over these rings.
\eproof

\subsection*{Functoriality}

The pairs $(G,A)$ where $G$ is a finite group that acts on a commutative ring $A$ can be made into a category such that a morphism 
\begin{equation}
\label{Eq:fGAHB}
f:(G,A)\to(H,B)
\end{equation}
consists of a group homomorphism $G\leftarrow H$ and an $H$-equivariant ring homomorphism $f:A\to B$; this is opposite to \cite[Ch.~III, \S8]{Brown:CohomologyofGroups} and \cite[6.7.6]{Weibel:Intro-HA}. A morphism of pairs $f:(G,A)\to(H,B)$ as in \eqref{Eq:fGAHB} induces a functor of tensor triangulated categories 
\[
f_*:D^-(AG)\to D^-(BH)
\]
or equivalently $f_*:K^-(AG\PProj)\to K^-(BH\PProj)$, the latter defined by
\begin{equation}
\label{f*KAGKBH}
f_*(P)=P\otimes_AB
\end{equation}
with diagonal $H$-action using the restriction under $G\leftarrow H$ on the first component. One verifies that the functor $f_*$ restricts to a functor
\begin{equation}
\label{Eq:f*-DAGAperf-DBHBperf}
f_*:D^b(AG)_{A\pperf}\to D^b(BH)_{B\pperf},
\end{equation}
which is given by \eqref{f*KAGKBH} on the objects of $D^b(AG)_{A\pproj}$. Under the isomorphism \eqref{Eq:End*DAGA-H*GA}, this functor gives a ring homomorphism
\begin{equation}
\label{Eq:f*HGAHHB}
f_*:H^*(G,A)\to H^*(H,B).
\end{equation}
This is the usual functoriality of group cohomology as a functor of two variables; see \cite[Ch.~III, \S8]{Brown:CohomologyofGroups} or \cite[6.7.6]{Weibel:Intro-HA}.

\bLe
\label{Le:DAGAperf-colim}
Let $G$ be a finite group and $A=\varinjlim_i A_i$ a filtered colimit of commutative rings with an action of $G$. Then 
\[
D^b(AG)_{A\pperf}\cong\varinjlim_iD^b(A_iG)_{A_i\pperf}
\] 
as a filtered $2$-colimit of tensor triangulated categories.
\eLe

\bproof
By Lemma \ref{Le:DAG-Aproj-Aperf} we can replace $D^b(AG)_{A\pperf}$ by $\EEE(A)=D^b(AG)_{A\pproj}$. The natural functor $\varinjlim_i\EEE(A_i)\to\EEE(A)$ is surjective on isomorphism classes because every complex in $\EEE(A)$ is determined by finite data of matrices over $A$ subject to finitely many relations. The functor is fully faithful because for $X,Y\in\EEE(A)$, after choosing a quasi-isomorphism $P\to X$ as in Lemma \ref{Le:DAGAperf-fin-res} we have $\Hom_{\EEE(A)}(X,Y)=\Hom_{K(AG)}(P,Y)$, which commutes with filtered colimits of $A$. 
\eproof

\subsection*{Finiteness conditions}

Let $A^G\subseteq A$ be the ring of $G$-invariant elements.

\bDe
\label{De:AG-noeth}
The pair $(G,A)$ will be called noetherian if the ring $A^G$ is noetherian and $A$ is a finite $A^G$-module.
\eDe

\bLe
\label{Le:AG-noeth}
If $A$ is an algebra of finite type over a noetherian subring $B$ of $A^G$, then $A^G$ is of finite type over $B$, and the pair $(G,A)$ is noetherian.
\eLe

\bproof
Since $A$ is integral over $A^G$ by \cite[Ch.\ V, \S 1.9, Prop.\ 22]{Bourbaki:AC}, \cite[Prop.\ 7.8]{AM:Commutative} implies that $A^G$ is a $B$-algebra of finite type. Hence $A^G$ is noetherian. Moreover $A$ is integral and of finite type over $A^G$, hence finite over $A^G$.
\eproof

\bCo
\label{Co:AG-ft-noeth}
If $A$ is a ring of finite type (i.e.\ a finitely generated $\ZZ$-algebra), then $A^G$ is of finite type as well, and the pair $(G,A)$ is noetherian.
\eCo

\bproof
Let $B$ be the image of $\ZZ\to A$ and apply Lemma \ref{Le:AG-noeth}.
\eproof

\bPr
\label{Pr:RAG-noeth}
If the pair $(G,A)$ is noetherian, the triangulated category $D^b(AG)_{A\pperf}$ is End-finite in the sense of Definition \ref{De:End-finite}, in particular the ring $H^*(G,A)$ is noetherian.
\ePr

\bproof 
Let $R=H^*(G,A^G)$. Since $A^G$ is a noetherian ring and $A$ is a finite $A^G$-module with an $A^G$-linear action of $G$, \cite[Thm.\ 6.1 \& Cor.~6.2]{Evens:Ring} give that $R$ is a noetherian ring and that $H^*(G,A)$ is a finite $R$-module, thus a noetherian ring. By Lemma \ref{Le:DAG-Aproj-Aperf} it suffices to show that for $P$, $Q \in D^b(AG)_{A\pproj}$ the graded $R$-module $M_{P,Q}=\Hom_{D(AG)}^*(P,Q)$ is finite. Since the finite $R$-modules form a Serre subcategory of the category of all $R$-modules, one can assume that the complexes $P$ and $Q$ are concentrated in degree zero. In that case, \eqref{Eq:End*DAGA-H*GA} gives $M_{P,Q}=H^*(G,\Hom_A(P,Q))$, which is finite over $R$ by \cite[Thm.\ 6.1]{Evens:Ring} again.
\eproof

\section{Perfect complexes over quotient stacks}

\label{Se:perfect}

As in \S \ref{Se:modules-skew} let $A$ be a commutative ring with an action of a finite group $G$. We consider the Deligne--Mumford stack $\XXX=[\Spec(A)/G]$ and the tensor triangulated category 
\begin{equation}
\label{Eq:TAG}
\TTT_{A,G}=\Perf(\XXX)
\end{equation}
of perfect complexes in $D(\XXX_{\et})$, the derived category of $\OOO_\XXX$-modules on the \'{e}tale site of $\XXX$, as well as the graded-commutative ring
\begin{equation}
\label{Eq:RAG}
R_{A,G}=\End^*_{\TTT_{A,G}}(\One)=\End^*_{D(\XXX_{\et})}(\OOO_\XXX).
\end{equation}
We refer to see \cite[(4.6.1)]{Laumon-MB:Champs} for the definition of $\XXX$, to \cite[\S 12]{Laumon-MB:Champs} for the \'{e}tale site of $\XXX$ and to \cite[\href{https://stacks.math.columbia.edu/tag/08G4}{Section 08G4}]{stacks-project} for perfect complexes on a ringed site. We note that $\Perf(\XXX)$ is the category of dualisable objects in $D(\XXX_{\et})$ by \cite[\href{https://stacks.math.columbia.edu/tag/0FPP}{Section 0FPP}]{stacks-project}. In particular, $\Perf(\XXX)$ is a rigid tensor triangulated category.

\bRe
For general algebraic stacks one uses the lisse-\'{e}tale site and the associated category $D(\XXX_{\liset})$ to define $\Perf(\XXX)$; see for example \cite{Hall-Rydh:Perfect}. For Deligne--Mumford stacks this makes no difference because the categories $D_{\qc}(\XXX_{\et})$ and $D_{\qc}(\XXX_{\liset})$ of complexes of $\OOO_\XXX$-modules with quasi-coherent cohomology are equivalent by \cite[Prop.~12.10.1]{Laumon-MB:Champs}; cf.\ also \cite[\S 1 and \S4]{Hall-Rydh:Perfect}.
\eRe

\bPr
\label{Pr:TAG-DAGAperf}
There is an equivalence of tensor triangulated categories 
\begin{equation}
\label{Eq:TAG-DAGAperf}
\TTT_{A,G}\cong D^b(AG)_{A\pperf}
\end{equation}
and hence an isomorphism of graded rings
\begin{equation}
\label{Eq:RAG-HGA}
R_{A,G}\cong H^*(G,A).
\end{equation}
\ePr

\bRe
\label{Re:TAG-DAGAperf}
We view Proposition \ref{Pr:TAG-DAGAperf} as a motivation to study the algebraic category $D^b(AG)_{A\pperf}$. In the remainder of the article, the stack $\XXX$ will not appear in an essential way, so the reader could skip the rest of this section and take \eqref{Eq:TAG-DAGAperf} as a definition.
\eRe

\bproof[Proof of Proposition \ref{Pr:TAG-DAGAperf}]
Let $Y=\Spec A$. We denote by $\OOO_Y[G]\MMod$ the category of $G$-equivariant $\OOO_Y$-modules on the \'{e}tale site of $Y$. Glueing of sheaves for the \'etale covering $\pi:Y\to\XXX$ yields an equivalence
\begin{equation}
\label{Eq:OX-OYG}
\OOO_\XXX\MMod\cong\OOO_Y[G]\MMod,
\end{equation}
using that $Y\times_\XXX Y\cong G\times Y$. There is a pair of adjoint functors
\begin{equation}
\label{Eq:varphi-Q}
\vcenter{
\xymatrix@M+0.2em{
AG\MMod \ar@<0.2em>[r]^-{\varphi} & \ar@<0.2em>[l]^-{Q} \OOO_Y[G]\MMod\cong\OOO_\XXX\MMod
}}
\end{equation}
with $\varphi$ left adjoint to $Q$, where $\varphi(M)=\tilde M$ is the quasi-coherent $\OOO_Y$-module associated to $M$ as an $A$-module, with the action of $G$ on $\tilde M$ induced by the given action on $M$, and $Q(\MMM)=\Gamma(Y,\MMM)$ as an $A$-module, carrying the action of $G$ induced by the given action on $\MMM$. The functor $\varphi$ is exact and preserves the tensor product defined by $M\otimes_AN$ in $AG\MMod$ and by $\MMM\otimes_{\OOO_Y}\NNN$ in $\OOO_Y[G]\MMod$, in both cases with diagonal $G$-action. We will show that $\varphi$ induces the inverse of the desired equivalence \eqref{Eq:TAG-DAGAperf}.

Let $D^+_{\qc}(\XXX_{\et})\subseteq D^+(\XXX_{\et})$ be the full subcategory of complexes with quasi-coherent cohomology. The following variant of \cite[II, Prop.~3.5]{SGA6} is a special case of \cite[Thm.~C.1]{Hall-Neeman-Rydh}; we give a direct proof for completeness.

\bLe
\label{Le:D+AG_D+qc}
The functor $\varphi$ induces an equivalence of triangulated categories $D^+(AG)\cong D^+_{\qc}(\XXX_{\et})$ with quasi-inverse functor $RQ$.
\eLe

\bproof
We begin with two initial remarks.

\medskip
(1)
For $M\in AG\MMod$,
the natural map $M\to Q(\varphi(M))$ is an isomorphism, so $\varphi$ is fully faithful. The image of $\varphi$ is the category of quasi-coherent $\OOO_\XXX$-modules in the sense of \cite[\href{https://stacks.math.columbia.edu/tag/03DL}{Definition 03DL}]{stacks-project} because an $\OOO_\XXX$-module $\MMM$ is quasi-coherent iff $\MMM_Y=\pi^*(\MMM)$ is a quasi-coherent $\OOO_Y$-module, which means that $\MMM_Y\cong \tilde M$ for an $A$-module $M$ by faithfully flat descent. 

\medskip
(2)
The functor $Q$ is left exact with derived functors $R^iQ(\MMM)=H^i(Y_{\et},\MMM)$ with the induced action of $G$. Since $Y$ is an affine scheme, for an $AG$-module $M$ and $i>0$ we have $H^i(Y_{\et},\tilde M)=0$, and thus $R^iQ(\varphi(M))=0$.

\medskip
Now the exact functor $\varphi$ induces an exact functor $\varphi:D^+(AG)\to D^+_{\qc}(\XXX_{\et})$. We have to show that for complexes $M\in D^+(AG)$ and $\MMM\in D^+_{\qc}(\XXX_{\et})$, the natural homomorphisms $\eta_M:M\to RQ(\varphi(M))$ and $\varepsilon_\MMM:\varphi(RQ(\MMM))\to\MMM$ are isomorphisms. Both assertions are easily reduced to the case where $M$ and $\MMM$ are modules concentrated in degree zero. Then $\eta_M$ is an isomorphism because $M\to Q(\varphi(M))$ is an isomorphism by (1) and $R^iQ(\varphi(M))=0$ for $i>0$ by (2). The module $\MMM$ is quasi-coherent since it lies in $D^+_{\qc}(\XXX_{\et})$, thus $\MMM=\varphi(N)$ for an $AG$-module $N$ by (1). Hence $\varepsilon_\MMM$ is an isomorphism since this holds for $\eta_N$.
\eproof

We continue the proof of Proposition \ref{Pr:TAG-DAGAperf}.
In the case $G=1$, by Lemma \ref{Le:D+AG_D+qc} the functor $\varphi$ gives an equivalence $D^+(A)\cong D^+_{\qc}(Y_{\et})$. This equivalence restricts to an equivalence $\Perf(A)\cong\Perf(Y_{\et})$ because a complex $P$ of $A$-modules is perfect iff for some faithfully flat ring homomorphism $A\to A'$ the complex $P\otimes_AA'$ is perfect; see \cite[\href{https://stacks.math.columbia.edu/tag/068T}{Lemma 068T}]{stacks-project}. For general $G$ it follows that a complex $P$ in $D^+(AG)$ is $A$-perfect iff the complex $\pi^*(\varphi(P))$ in $D^+(Y_{\et})$ is perfect, which means that $\varphi(P)$ in $D^+(\XXX_{\et})$ is perfect because this is a local condition. Hence the equivalence of Lemma \ref{Le:D+AG_D+qc} restricts to the desired equivalence $D^b(AG)_{A\pperf}\cong\Perf(\XXX)$. This functor preserves the tensor product because this evidently holds for its restriction to $D^b(AG)_{A\pproj}$; see Lemma \ref{Le:DAG-Aproj-Aperf}. The isomorphism \eqref{Eq:RAG-HGA} follows by \eqref{Eq:End*DAGA-H*GA}.
\end{proof}

\subsection*{Functoriality}

A morphism of pairs $f:(G,A)\to(H,B)$ as in \eqref{Eq:fGAHB} induces a morphism of algebraic stacks (see for example \cite[\href{https://stacks.math.columbia.edu/tag/046Q}{Lemma 046Q}]{stacks-project})
\[
\psi=\Spec(f):[\Spec(B)/H]\to[\Spec(A)/G],
\]
which gives an inverse image functor of tensor triangulated categories
\begin{equation}
\label{Eq:phi*TAG-TBH}
f_*=\psi^*:\TTT_{A,G}\to\TTT_{B,H},
\end{equation}
using that the \'etale topos of Deligne--Mumford stacks is functorial, see for example \cite[Constr.~2.4]{Zheng:SixOpDM}, together with \cite[\href{https://stacks.math.columbia.edu/tag/08H6}{Lemma 08H6}]{stacks-project}. One verifies that under the equivalence of Proposition \ref{Pr:TAG-DAGAperf}, the functor $f_*$ of \eqref{Eq:phi*TAG-TBH} corresponds to the functor $f_*$ of \eqref{Eq:f*-DAGAperf-DBHBperf}. We leave out further details.

\section{The comparison map: basic properties}
\label{Se:comparison-basic}

As earlier let $A$ be a commutative ring with an action of a finite group $G$. The comparison map $\rho_\TTT$ of \eqref{Eq:rho} for the category $\TTT=\TTT_{A,G}$ of \eqref{Eq:TAG} will be denoted by
\begin{equation}
\label{Eq:rhoAG}
\rho_{A,G}:\Spc(\TTT_{A,G})\to\Spec^h(R_{A,G}).
\end{equation}

\subsection*{Functoriality} 

Since the comparison map $\rho_\TTT$ is natural in $\TTT$, for a morphism of pairs $f:(G,A)\to(H,B)$ as in \eqref{Eq:fGAHB} the functor $f_*:\TTT_{A,G}\to\TTT_{B,H}$ of \eqref{Eq:phi*TAG-TBH} induces a commutative diagram of topological spaces
\begin{equation}
\label{Dia:TBH-TAG}
\vcenter{
\xymatrix@M+0.2em@C+1em{
\Spc(\TTT_{B,H}) \ar[r]^-{\rho_{B,H}} \ar[d]_{f^\TTT} & \Spec^h(R_{B,H}) \ar[d]^{f^R} \\
\Spc(\TTT_{A,G}) \ar[r]^-{\rho_{A,G}} & \Spec^h(R_{A,G}).
}}
\end{equation}
Here $f^\TTT$ is the inverse image map under the functor $f_*$, and $f^R$ is the inverse image map under the ring homomorphism $f_*:R_{A,G}\to R_{B,H}$ defined by this functor. We recall that under the equivalence \eqref{Eq:TAG-DAGAperf}, the functor $f_*$ corresponds to the functor \eqref{Eq:f*-DAGAperf-DBHBperf}, and under the isomorphism \eqref{Eq:RAG-HGA}, the ring homomorphism $f_*$ corresponds to the homomorphism \eqref{Eq:f*HGAHHB}.

\subsection*{The case of the trivial group}

We will write $\TTT_A=\TTT_{A,\{e\}}=\Perf(\Spec A)$ and $\rho_A=\rho_{A,\{e\}}$. In this case we have $R_{A,\{e\}}=A$ as a graded ring concentrated in degree zero, and the comparison map
\begin{equation}
\label{Eq:rhoA}
\rho_A:\Spc(\TTT_A)\to\Spec(A)
\end{equation}
is a homeomorphism by \cite[Prop.\ 8.1]{Balmer:Spectra3}, which is based on \cite[Thm.\ 3.15]{Thomason:Classification}.

\subsection*{Restriction to fibers}

For a general pair $(G,A)$ there is a morphism of pairs $\pi:(\{e\},A^G)\to(G,A)$ defined by the inclusion $A^G\to A$ and the unique group homomorphism $\{e\}\leftarrow G$. This gives the following instance of the functoriality diagram \eqref{Dia:TBH-TAG}.
\begin{equation}
\label{Dia:fiber}
\vcenter{
\xymatrix@M+0.2em@C+2em{
\Spc(\TTT_{A,G}) \ar[r]^{\rho_{A,G}} \ar[d]_{\pi^\TTT} &  \Spec^h(R_{A,G}) \ar[d]^{\pi^R} \\
\Spc(\TTT_{A^G}) \ar[r]_{\rho_{A^G}}^\sim & \Spec (A^G)
}}
\end{equation}
Here $\rho_{A^G}$ is a homeomorphism by \cite[Prop.\ 8.1]{Balmer:Spectra3} as explained in \eqref{Eq:rhoA}. We note that $\pi^R$ is induced by the ring homomorphism $\pi_*:A^G\to R_{A,G}$ given by the inclusion of the degree zero component, and hence $\pi^R(\Fp)=\Fp^0$ is the degree zero part of a homogeneous prime ideal $\Fp$ of $R_{A,G}$.

For a given $\Fq\in\Spec(A^G)$ with unique inverse image $\QQQ\in\Spc(\TTT_{A^G})$ the vertical fibers in \eqref{Dia:fiber} over these points will be denoted by
\begin{equation}
\label{Eq:fibers}
\Spc(\TTT_{A,G})_\Fq=(\pi^\TTT)^{-1}(\QQQ)
\quad\text{and}\quad
\Spec^h(R_{A,G})_\Fq=(\pi^R)^{-1}(\Fq).
\end{equation}
The map $\rho_{A,G}$ induces a map between these fibers
\begin{equation}
\label{Eq:rhoAGq}
(\rho_{A,G})_\Fq:\Spc(\TTT_{A,G})_\Fq\to\Spec^h(R_{A,G})_\Fq.
\end{equation}
The following is evident.

\bLe
\label{Le:rhobij-rhoqbij}
For a given pair $(G,A)$, the map $\rho_{A,G}$ is bijective iff the map $(\rho_{A,G})_\Fq$ is bijective for each $\Fq\in\Spec(A^G)$.
\qed
\eLe

\bRe
The diagram \eqref{Dia:fiber} is functorial with respect to $(G,A)$. More precisely, a morphism of pairs $f:(G,A)\to(H,B)$ as in \eqref{Eq:fGAHB} gives rise to a commutative cube: $[\eqref{Dia:fiber} \text{ for } (H,B)] \to [\eqref{Dia:fiber} \text{ for } (G,A)]$. Since the lower line of \eqref{Dia:fiber} is always a homeomorphism, the essential information of this cube is captured by the following extension of \eqref{Dia:TBH-TAG},
\begin{equation}
\label{Dia:TBH-TAG-extended}
\vcenter{
\xymatrix@M+0.2em{
\Spc(\TTT_{B,H})
\ar[d]_{f^\TTT} \ar[r]^-{\rho_{B,H}} &
\Spec^h(R_{B,H})
\ar[d]^{f^R} \ar[r]^-{\pi^R} &
\Spec(B^H) 
\ar[d]^{f^0} \\
\Spc(\TTT_{A,G})
\ar[r]^-{\rho_{A,G}} &
\Spec^h(R_{A,G})
\ar[r]^-{\pi^R} &
\Spec(A^G),
}}
\end{equation}
where $f^0$ comes from the ring homomorphism $A^G\to B^H$ induced by $f$.
\eRe

\subsection*{Fibers of the coefficient ring}

For a given $\Fq\in\Spec(A^G)$ let 
\begin{equation}
\label{Eq:A(q)}
A(\Fq)=(A\otimes_{A^G}k(\Fq))_{\red}
\end{equation}
where the index ${\red}$ means maximal reduced quotient, so $\Spec A(\Fq)$ is the reduced fiber over $\Fq$ of the morphism $\Spec A\to\Spec A^G$, and let 
\begin{equation}
\label{Eq:psi_q}
\psi_\Fq:A\to A(\Fq)
\end{equation}
be the natural homomorphism given by $a\mapsto a\otimes 1$. The action of $G$ on $A$ induces an action on $A(\Fq)$ via the first factor, and $\psi_\Fq$ is $G$-equivariant.

\bLe
\label{Le:Aqkpi}
Let $\Fp_1,\ldots,\Fp_r$ be the prime ideals of $A$ lying over $\Fq$. Then
\[
A(\Fq)=k(\Fp_1)\times\ldots\times k(\Fp_r).
\]
\eLe

\bproof
The prime ideals of $A$ over $\Fq$ form a finite discrete set because they form a single $G$-orbit in $\Spec(A)$ by \cite[Ch.\ V, \S 2.2, Thm.\ 2 (i)]{Bourbaki:AC}, and this set is homeomorphic to $\Spec(A(\Fq))$. Hence the reduced ring $A(\Fq)$ is the product of its residue fields, and these residue fields coincide with the corresponding residue fields of $A$ since $A(\Fq)$ is a quotient of a localisation of $A$.
\eproof

\bLe
\label{Le:AG-AqG-factors}
There is a commutative diagram of rings
\begin{equation*}
\vcenter{
\xymatrix@M+0.2em{
A \ar[rr]^{\psi_{\Fq}} &&
A(\Fq) \\
A^G \ar[u] \ar[r]^{can} &
k(\Fq) \ar[r] & A(\Fq)^G \ar[u] 
}}
\end{equation*}
where the vertical homomorphisms are the inclusions. Here $A(\Fq)^G$ is a field, and $k(\Fq)\to A(\Fq)^G$ is a purely inseparable field extension.
\eLe

\bproof
The composition $A^G\to A\to A(\Fq)$ factors over $k(\Fq)$ by the definition of $A(\Fq)$, and the resulting homomorphism $k(\Fq)\to A(\Fq)$ has image in $A(\Fq)^G$ since $G$ acts trivially on $k(\Fq)$. If $\Fp\in\Spec(A)$ lies over $\Fq$ and if $H\subseteq G$ is the stabiliser of $\Fp$, then $A(\Fq)^G\cong k(\Fp)^H$, which is a field, and a purely inseparable extension of $k(\Fq)$ by \cite[Ch.\ V, \S 2.2, Thm.\ 2 (ii)]{Bourbaki:AC}.
\eproof

\subsection*{The fiber diagram}

Again let $\Fq\in\Spec(A^G)$ be given.

\bPr
The functoriality diagram \eqref{Dia:TBH-TAG} for the homomorphism of pairs $\psi_\Fq:(G,A)\to(G,A(\Fq))$ of \eqref{Eq:psi_q} induces a commutative diagram
\begin{equation}
\label{Dia:fiber-diagram}
\vcenter{
\xymatrix@M+0.2em@C+1em{
\Spc(\TTT_{A(\Fq),G}) \ar[r]^-{\rho_{A(\Fq),G}} \ar[d]_{\psi_{\Fq,\res}^\TTT} & \Spec^h(R_{A(\Fq),G}) \ar[d]^{\psi_{\Fq,\res}^R} \\
\Spc(\TTT_{A,G})_\Fq \ar[r]^-{(\rho_{A,G})_\Fq} & \Spec^h(R_{A,G})_\Fq,
}}
\end{equation}
which we call the fiber diagram at $\Fq$.
\ePr

This diagram also appears in the introduction, where the vertical arrows are denoted by $j_\TTT$ and $j_R$.

\bproof
Let us draw the extended functoriality diagram \eqref{Dia:TBH-TAG-extended} for $\psi_q$.
\begin{equation*}
\vcenter{
\xymatrix@M+0.2em@C+1em{
\Spc(\TTT_{A(\Fq),G}) \ar[r]^-{\rho_{A(\Fq),G}} \ar[d]_{\psi_\Fq^{\TTT}} & \Spec^h(R_{A(\Fq),G}) \ar[d]^{\psi_\Fq^{R}} \ar[r]^-{\pi_R} &
\Spec(A(\Fq)^G) \ar[d]^{\psi_\Fq^0} \\
\Spc(\TTT_{A,G}) \ar[r]^-{\rho_{A,G}} & 
\Spec^h(R_{A,G}) \ar[r]^-{\pi_R} &
\Spec(A^G)
}}
\end{equation*}
We have to verify that the images of the vertical arrows $\psi_\Fq^\TTT$ and $\psi_\Fq^R$ map to $\Fq$ in $\Spec(A^G)$. The holds because the ring homomorphism $A^G\to A(\Fq)^G$ factors over $k(\Fq)$ by Lemma \ref{Le:AG-AqG-factors}, so the image of $\psi_\Fq^0$ is the singleton $\{\Fq\}$.
\eproof

\bLe
\label{Le:TAqG-TLH}
As earlier let $\Fp_1,\ldots,\Fp_r$ be the prime ideals of $A$ over $\Fq$.
Moreover let $L=k(\Fp_1)$, and let $H\subseteq G$ be the stabiliser of the element $\Fp_1$ of $\Spec(A)$. Then the tensor triangulated category $\TTT_{A(\Fq),G}$ is equivalent to $\TTT_{L,H}$, and consequently the map $\rho_{A(\Fq),G}$ is isomorphic to $\rho_{L,H}$.
\eLe

\bproof
The group $H$ acts on $L$ by functoriality. We use Lemma \ref{Le:Aqkpi}. Since $G$ acts transitively on the set $\{\Fp_1,\ldots,\Fp_r\}$, there is an equivalence 
\[
A(\Fq)G\mmod\cong LH\mmod
\] 
given by $M\mapsto M\otimes_{A(\Fq)}L$ with diagonal action of $H$. The resulting equivalence $D^b(A(\Fq)G\mmod)\cong D^b(LH\mmod)$ gives $\TTT_{A(\Fq),G}\cong\TTT_{L,H}$ by Proposition \ref{Pr:TAG-DAGAperf}, using that every $A(\Fq)$-module is projective.
\eproof

\bRe
On the geometric side, the equivalence $\TTT_{A(\Fq),G}\cong\TTT_{L,H}$ comes from an isomorphism of stacks $[\Spec(L)/H]\cong[\Spec(A(\Fq))/G]$ which is induced by the obvious morphism of pairs $(G,A(\Fq))\to (H,L)$.
\eRe

\subsection*{Additional comments on the ring $A(\Fq)$}

\bLe
\label{Le:AG-local}
If the ring $A^G$ is local with maximal ideal $\Fq$, then $\psi_\Fq$ is surjective and induces a homeomorphism $\Spec(A(\Fq))\cong\Max(A)$.
\eLe

\bproof
The assumption implies that $A^G\to k(\Fq)$ is surjective, and hence $\psi_\Fq$ is surjective. Since $A^G\subseteq A$ is an integral extension, a prime ideal $\Fp$ of $A$ is maximal iff $\Fp\cap A^G$ is a maximal ideal of $A^G$; see \cite[Ch.~V, \S2.1, Prop.~1]{Bourbaki:AC}. Hence $\Max(A)$ is the set of prime ideals of $A$ lying over $\Fq$, which is homeomorphic to $\Spec(A(\Fq))$. 
\eproof

\bLe
\label{Le:AqBq-surj}
Let $f:A\to B$ be a $G$-equivariant ring homomorphism and let $\tilde\Fq\in\Spec(B^G)$ with image $\Fq\in\Spec(A^G)$ be given. Then $f$ induces an injective ring homomorphism $f':A(\Fq)\to B(\tilde\Fq)$. If $B$ is a localisation of a quotient of $A$, then $f'$ is bijective. 
\eLe

\bproof
Clearly $f$ induces $f'$. All assertions follow from Lemma \ref{Le:Aqkpi}. Indeed, the $G$-equivariant map $\Spec(f):\Spec(B)\to\Spec(A)$ sends the $G$-orbit over $\tilde\Fq$ to the $G$-orbit over $\Fq$, and this map between $G$-orbits is necessarily surjective; moreover for $\tilde\Fp\in\Spec(B)$ with image $\Fp\in\Spec(A)$ the homomorphism of residue fields $k(\Fp)\to k(\tilde\Fp)$ is injective; hence $f'$ is injective. If $B$ is a localisation of a quotient of $A$, then $\Spec(f)$ is injective, so our map between $G$-orbits is bijective, moreover $f$ induces isomorphisms of the residue fields; hence $f'$ is bijective.
\eproof

\bLe
\label{Le:psiq-factors}
Let $A_\Fq=S^{-1}A$ with $S=A^G\setminus\Fq$. The homomorphism $\psi_\Fq$ factors into $G$-equivariant homomorphisms
\[
A\to A_\Fq\xrightarrow\pi A(\Fq)
\]
where $\pi$ is surjective. The ring $(A_{\Fq})^G=(A^G)_{\Fq}$ is local with maximal ideal $\Fq_{\Fq}$. There is an isomorphism $A(\Fq)\cong A_\Fq(\Fq_\Fq)$ under which $\pi$ corresponds to the homomorphism \eqref{Eq:psi_q} for $A_\Fq$ and $\Fq_\Fq$ in place of $A$ and $\Fq$.
\eLe

\bproof
The homomorphism $\psi_\Fq$ factors over $S^{-1}A$ because $A^G\to k(\Fq)$ factors over $S^{-1}(A^G)$. We have $(A_{\Fq})^G=(A^G)_{\Fq}$ since localisation is flat. The rest follows easily; one can also use Lemma \ref{Le:AqBq-surj} with $B=S^{-1}A$.
\eproof

\section{Spectra of graded rings}

For a graded-commutative ring $R$ let $R_{\even}$ be the subring of $R$ generated by the homogeneous elements of even degree. Then $R_{\even}$ is a commutative graded ring, and there is a homeomorphism $\Spec^h(R)\cong\Spec^h(R_{\even})$. In this section we record a number of basic properties of this construction.

\bLe
\label{Le:fh-surj-homeo}
For a homomorphism $S\to T$ of graded-commutative rings let
\[
f:\Spec(T_{\even})\to\Spec(S_{\even})\quad \text{and}\quad f_h:\Spec^h(T)\to\Spec^h(S)
\]
be the induced maps. If $f$ is surjective, then so is $f_h$. If $f$ is a homeomorphism, then so is $f_h$.
\eLe

\bproof
One can replace $S$ and $T$ by $S_{\even}$ and $T_{\even}$. Since $\Spec^h(S)\subseteq\Spec(S)$ and $\Spec^h(T)\subseteq\Spec(T)$ carry the subspace topology, it suffices to show that $f_h$ is surjective if $f$ is surjective. For $\Fp\in\Spec^h(S)$ we consider the residue field $k(\Fp)=\Frac(S/\Fp)$ and the graded residue field $k((\Fp))=S_{\Fp}/\Fp S_{\Fp}$, where $S_\Fp$ is the graded localisation of $S$ at $\Fp$. There are ring homomorphisms $S\to k((\Fp))\to k(\Fp)$. Since $f$ is surjective, the ring $T\otimes_Sk(\Fp)$ is non-zero, hence the graded ring $T\otimes_Sk((\Fp))$ is non-zero. Any graded prime ideal of $T\otimes_Sk((\Fp))$ gives a graded prime ideal of $T$ that maps to $\Fp$.
\eproof

\bLe
\label{Le:graded-invariants}
Let $R$ be a graded-commutative ring with an action of a finite group $\Gamma$ by automorphisms of graded rings. Then the inclusion $R^\Gamma\to R$ induces a homeomorphism $\Spec^h(R)/\Gamma\cong\Spec^h(R^\Gamma)$.
\eLe

\bproof
The lemma includes its non-graded version because a commutative ring can be considered as a graded ring concentrated in degree zero. The non-graded version is well-known and appears for example in \cite[Exp.~V, Prop.~1.1]{SGA1}. In more detail, $R^\Gamma\to R$ is integral and gives a bijective continuous map $\Spec(R)/\Gamma\to\Spec(R^\Gamma)$ by \cite[Ch.~V, \S 2, Th.~2]{Bourbaki:AC}; this map is also closed by \cite[Prop.~5.12]{GW:AG1}. In the graded case it follows that the natural map $\Spec^h(R)/\Gamma\to\Spec^h(R^\Gamma)$ is the inclusion of a subspace, and the map is surjective by Lemma \ref{Le:fh-surj-homeo}, hence a homeomorphism.
\eproof

\bLe
\label{Le:graded-noetherian}
A graded commutative ring $R=\bigoplus_{n\ge 0}R_n$ is noetherian iff the ring $R_0$ is noetherian and $R$ is an $R_0$-algebra of finite type. If this holds, for any $d>0$ the ring $R^{(d)}=\bigoplus_nR_{nd}$ is noetherian as well, and $R$ is a finite $R^{(d)}$-module.
\eLe

\bproof
See for example \cite[Thm.~13.1]{Matsumura:CRT} for the first assertion. The second assertion is reduced to the case $R=R_0[T_1,\ldots,T_r]$ where each $T_i$ is homogeneous of some positive degree. Then $R'=R_0[T_1^d,\ldots,T_r^d]\subseteq R^{(d)}\subseteq R$ where $R'$ is noetherian and $R$ is finite over $R'$, and the second assertion follows. \eproof

\section{The comparison map: the field case}
\label{Se:field}

Let $L$ be a field with an action of a finite group $G$. By Proposition \ref{Pr:TAG-DAGAperf}, the tensor triangulated category $\TTT_{L,G}$ is equivalent to $D^b(LG\mmod)$. The aim of this section is to prove the following result.

\begin{Thm}
\label{Th:field}
If $L$ is a field, the map $\rho_{L,G}:\Spc(\TTT_{L,G})\to\Spec^h(R_{L,G})$ of \eqref{Eq:rhoAG} is a homeomorphism.
\end{Thm}

If $G$ acts trivally on $L$, this is the content of \cite[Prop.\ 8.5]{Balmer:Spectra3}, which is eventually based on the fact that thick tensor ideals in $D^b(LG\mmod)$ are classified by their cohomological support in $\Spec^h(R_{L,G})$, or equivalently that  thick tensor ideals in the stable module category of $G$ over $L$ are classified by their cohomological support in $\Proj(R_{L,G})$. The stable version is proved in \cite[Thm.\ 3.4]{BCR:Thick} when $L$ is algebraically closed and in \cite[Thm.~11.4]{BIK:Stratifying} for arbitrary fields $L$; a direct proof of the unstable version appears in \cite{Carlson-Iyengar:Thick}.

The general case will be reduced to the case of trivial action. Let 
\[
H=\Ker(G\to\Aut(L))
\]
so that $G/H$ acts faithfully on $L$. The morphism of pairs 
\[
f:(G,L)\to(H,L)
\] 
defined by the identity of $L$ and the inclusion $G\leftarrow H$ gives the following instance of the functoriality diagram \eqref{Dia:TBH-TAG}.
\begin{equation}
\label{Dia:TLH-TLG}
\vcenter{
\xymatrix@M+0.2em@C+1em{
\Spc(\TTT_{L,H}) \ar[r]^-{\rho_{L,H}} \ar[d]_{f^\TTT} & \Spec^h(R_{L,H}) \ar[d]^{f^R} \\
\Spc(\TTT_{L,G}) \ar[r]^-{\rho_{L,G}} & \Spec^h(R_{L,G})
}}
\end{equation}
Here the vertical arrows are induced by the functor $f_*:\TTT_{L,G}\to\TTT_{L,H}$ which corresponds to the restriction functor
\[
\res_H^G:D^b(LG\mmod)\to D^b(LH\mmod)
\]
under the equivalence of Proposition \ref{Pr:TAG-DAGAperf}.

\bLe
\label{Le:field1}
The group $\bar G=G/H$ acts on all spaces in \eqref{Dia:TLH-TLG} with trivial action on the lower line such that all maps in \eqref{Dia:TLH-TLG} are $\bar G$-equivariant.
\eLe

\bproof
The actions are induced by the conjugation action of $G$ on the ring $LG$ and on the subring $LH$; note that $H$ is a normal subgroup of $G$.

In more detail, for $g\in G$ and an $LH$-module $X$ we form the $LH$-module $X^g$ which is $X$ with the action of $z\in LH$ by $gzg^{-1}$. This defines a right action of $G$ on the triangulated category $\TTT_{L,H}\cong D^b(LH\mmod)$, called here the conjugation action. The conjugation action admits the following alternative description, which shows that the action of $G$ on $\TTT_{L,H}$ preserves the tensor structure: Each $g\in G$ gives a homomorphism of pairs $(H,L)\to(H,L)$ defined by $L\to L$, $a\mapsto g^{-1}(a)$ and $H\leftarrow H$, $ghg^{-1}\mapsfrom h$. The resulting endomorphism of $\TTT_{L,H}$ by functoriality with respect to $(H,L)$ is isomorphic to the endomorphism $X\mapsto X^g$ because there is an isomorphism $X^g\to X\otimes_{L,g^{-1}}L$, $y\mapsto y\otimes 1$. Hence the action of $G$ on the pair $(H,L)$ induces by functoriality the conjugation action of $G$ on $\TTT_{L,H}$. 

The conjugation action of $G$ on $\TTT_{L,H}$ induces compatible left actions of $G$ on the source and target of $\rho_{L,H}$.
Similarly, $G$ acts on the tensor triangulated category $\TTT_{L,G}$, which induces compatible actions of $G$ on the source and target of $\rho_{L,G}$ such that all maps in \eqref{Dia:TLH-TLG} are equivariant. 
It remains to verify that $H$ acts trivially on the source and target of $\rho_{L,H}$ and $G$ acts trivially on the source and target of $\rho_{L,G}$.

For an $LG$-module $X$ and $g\in G$, multiplication by $g$ is an isomorphism $X\cong X^g$. It follows that $G$ acts trivially on $\Spc(\TTT_{L,G})$. The alternative description of the conjugation action implies that the resulting action of $G$ on $R_{L,G}\cong H^*(G,L)$ corresponds to the conjugation action in group cohomology as defined in \cite[Ch.~III, \S8]{Brown:CohomologyofGroups}, which is trivial by \cite[Ch.~III, Prop.~(8.1)]{Brown:CohomologyofGroups}. Similarly, $H$ acts trivially on $\Spc(\TTT_{L,H})$ and on $R_{L,H}$.
\eproof

\bLe
\label{Le:field2}
The map $f^R$ of \eqref{Dia:TLH-TLG} induces a homeomorphism
\[
\Spec^h(R_{L,H})/\bar G\xrightarrow{\;\;\bar f^R\;\;}\Spec^h(R_{L,G}).
\]
\eLe

\bproof
Let $K=L^G$. Then $L/K$ is a finite Galois extension with Galois group $\bar G$. We consider the sequence of left exact functors
\begin{equation}
\label{Eq:Functors-HS}
LG\MMod\xrightarrow{(-)^H}L\bar G\MMod\xrightarrow{(-)^{\bar G}}K\MMod
\end{equation}
whose composition is the functor of $G$-invariants. By Galois descent \cite[Thm.~14.85]{GW:AG1}, the functor $(-)^{\bar G}$ in \eqref{Eq:Functors-HS} is an equivalence, in particular exact. Hence \eqref{Eq:Functors-HS} yields an isomorphism of $\delta$-functors
\[
H^i(G,M)\xrightarrow\sim H^i(H,M)^{\bar G}
\]
for $M\in LG\MMod$, which coincides with the restriction map in group cohomology because this holds for $i=0$ as is easily verified. For $M=L$ it follows that the homomorphism $f_*:R_{L,G}\to R_{L,H}$ identifies $R_{L,G}$ with the ring of $\bar G$-invariants in $R_{L,H}$. Then Lemma \ref{Le:graded-invariants} finishes the proof.
\eproof

\bLe
\label{Le:field3}
The map $f^\TTT$ of \eqref{Dia:TLH-TLG} is surjective.
\eLe

\bproof
The restriction functor $f_*:LG\mmod\to LH\mmod$ has an exact right adjoint $f^!$ defined by $f^!(M)=\Hom_{LH}(LG,M)$, which is an $LG$-module using the right $LG$-module structure of $LG$. This induces an exact right adjoint $f^!$ of the tensor triangulated functor $f_*:D^b(LG\mmod)\to D^b(LH\mmod)$. By \cite[Thm.~1.7]{Balmer:Surjectivity} it follows that the image of the map $f^\TTT=\Spc(f_*)$ is equal to the support of $f^!(\One)=\Hom_{LH}(LG,L)$, viewed as an object of $D^b(LG\mmod)$. One verifies that $H$ acts trivially on $f^!(\One)$, so this is an $L\bar G$-module. An $L\bar G$-module is determined by its $L$-dimension by Galois descent \cite[Thm.~14.85]{GW:AG1}. Hence $f^!(\One)\cong L^{[G:H]}$ and thus $\supp(f^!(\One))=\supp(\One)=\Spc(\TTT_{L,G})$.
\eproof

\bproof[Proof of Theorem \ref{Th:field}]
By Lemma \ref{Le:field1}, the commutative diagram \eqref{Dia:TLH-TLG} induces a commutative diagram of topological spaces
\begin{equation}
\label{Dia:TLH-TLG-Quot}
\vcenter{
\xymatrix@M+0.2em@C+1em{
\Spc(\TTT_{L,H})/\bar G \ar[r]^-{\tilde \rho_{L,H}} \ar[d]_{\tilde f^\TTT} & \Spec^h(R_{L,H})/\bar G \ar[d]^{\tilde f^R} \\
\Spc(\TTT_{L,G}) \ar[r]^-{\rho_{L,G}} & \Spec^h(R_{L,G}).
}}
\end{equation}
Here $\tilde \rho_{L,H}$ is a homeomorphism because $\rho_{L,H}$ is a homeomorphism by \cite[Prop.\ 8.5]{Balmer:Spectra3}, $\tilde f^R$ is a homeomorphism by Lemma \ref{Le:field2}, and $\tilde f^\TTT$ is surjective by Lemma \ref{Le:field3}. It follows that all arrows in \eqref{Dia:TLH-TLG-Quot} are homeomorphisms.
\eproof

\bRe
\label{Re:field}
In the situation of Theorem \ref{Th:field} one could ask for a relation between the stacks $[\Spec(L)/G]$ and $[\Spec(K)/H]$ or between the corresponding module categories $LG\mmod$ and $KH\mmod$. In general, these stacks are not isomorphic, and the categories are not equivalent.
Example: $L=\CC$ and $G=\ZZ/4\ZZ$ such that a generator of $G$ acts on $L$ by complex conjugation, so $K=\RR$ and $H=2\ZZ/4\ZZ$. Then $KH\mmod$ has two isomorphism classes of one-dimensional representations while $LG\mmod$ has only one such class, so these categories are not equivalent as tensor categories because the number of invertible objects up to isomorphism is different.
\eRe

\section{Change of coefficients in group cohomology}
\label{Se:change-coefficient}

We fix a finite group $G$. For a commutative ring $A$ with an action of $G$ we ask how the graded-commutative ring $R_{A,G}=H^*(G,A)$ and its homogeneous prime spectrum depend on $A$.

\subsection*{Base change homomorphisms}

A $G$-equivariant ring homomorphism
\[
u:A\to B
\]
induces a homomorphism of graded-commutative rings
\begin{equation*}
u':R_{A,G}\to R_{B,G}
\end{equation*}
whose degree zero component is the homomorphism $A^G\to B^G$ defined by $u$. This gives a homomorphism of graded-commutative rings
\begin{equation*}
u'':R_{A,G}\otimes_{A^G}B^G\to R_{B,G}
\end{equation*}
and, by restriction to the subrings generated by the elements of even degree, a homomorphism of commutative graded rings
\begin{equation}
\label{Eq:RAGevRBGev}
u''':(R_{A,G})_{\even}\otimes_{A^G}B^G\to(R_{B,G})_{\even}.
\end{equation}

A homomorphism of commutative rings $S\to T$ will be called a universal homeomorphism if it indues a universal homeomorphism $\Spec T\to\Spec S$ of schemes, which means that for every ring homomorphism $S\to S'$ the natural map $\Spec(T\otimes_SS')\to\Spec(S')$ is a homeomorphism. 

\bDe
\label{De:Coh-uh}
We denote by $\Cohuh(G)$ the class of all $G$-equivariant homomorphisms of commutative rings $u:A\to B$ such that the base change homomorphism $u'''$ in \eqref{Eq:RAGevRBGev} is a universal homeomorphism.
\eDe

\bRe
\label{Re:Coh-uh}
The class $\Cohuh(G)$ does not contain all $G$-equivariant ring homomorphisms. For example, if $L/K$ is a finite Galois extension with Galois group $G$, then $L$ is an induced $G$-module by the normal basis theorem, the inclusion map $K\to L$ is $G$-equivariant, and the associated homomorphism $R_{K,G}\to R_{L,G}$ is the augmentation $H^*(G,K)\to K$.
\eRe

\bPr
\label{Pr:Coh-uh-comp}
The class $\Cohuh(G)$ is stable under composition.
\ePr

\bproof
If $A\xrightarrow u B\xrightarrow v C$ is a sequence of $G$-equivariant ring homomorphisms, then $(v\circ u)'''=v'''\circ(u'''\otimes_{B^G}C^G)$. \eproof

We will show that the class $\Cohuh(G)$ contains all $G$-equivariant surjections and localisations. First we study the effect of such homomorphisms on the rings of $G$-invariants, i.e.\ on the degree zero part of $R_{A,G}$.

\bPr
\label{Pr:AGIGBG}
Let $A$ be a commutative ring with an action of $G$ and let $B=A/I$ for a $G$-stable ideal $I$ of $A$. Then the natural ring homomorphism $A^G/I^G\to B^G$ is a universal homeomorphism.
\ePr

\bproof
Since $A$ is integral over $A^G$ by \cite[Ch.\ V, \S 1.9, Prop.\ 22]{Bourbaki:AC}, the injective ring homomorphism $A^G/I^G\to B$ is integral, so $A^G/I^G\to B^G$ is integral and injective as well, and the induced morphism
\[
f:\Spec(B^G)\to\Spec(A^G/I^G)
\] 
is integral and surjective by \cite[\href{https://stacks.math.columbia.edu/tag/00GQ}{Lemma 00GQ}]{stacks-project}. We have
\begin{equation*}
\Spec(A^G)=\Spec(A)/G
\qquad
\text{and}
\qquad
\Spec(B^G)=\Spec(B)/G
\end{equation*}
as topological spaces; see Lemma \ref{Le:graded-invariants}. Since $\Spec(B)\to\Spec(A)$ is injective it follows that $\Spec(B^G)\to\Spec(A^G)$ is injective, hence $f$ is injective and thus bijective. Since $f$ is integral and bijective, by \cite[18.12.11]{EGAIV4} it remains to show that the residue field extensions induced by $f$ are purely inseparable.

For given $\tilde\Fq\in\Spec(B^G)$ with image $\Fq\in\Spec(A^G)$ the homomorphism $A\to B$ induces an isomorphism $A(\Fq)\to B(\tilde\Fq)$ by Lemma \ref{Le:AqBq-surj}. The sequence of $G$-equivariant ring homomorphisms 
\[
A\to B\to B(\tilde\Fq)\cong A(\Fq)
\]
gives ring homomorphisms $A^G\to B^G\to B(\tilde\Fq)^G\cong A(\Fq)^G$, which induces field extensions $k(\Fq)\to k(\tilde\Fq)\to B(\tilde\Fq)^G\cong A(\Fq)^G$ by Lemma \ref{Le:AG-AqG-factors} applied to $B$, and the total extension $k(\Fq)\to A(\Fq)^G$ is purely inseparable by Lemma \ref{Le:AG-AqG-factors} applied to $A$. Hence $k(\Fq)\to k(\tilde\Fq)$ is purely inseparable as well.
\eproof

\bLe
\label{Le:AGBGloc}
Let $A\to B$ be a $G$-equivariant homomorphism of rings such that $B$ is a localisation of $A$, and let $S\subseteq A$ be the set of elements which become invertible in $B$. Then $B^G=(S^G)^{-1}A^G$ and $B=A\otimes_{A^G}B^G$.
\eLe

\bproof
If $s\in S$ then $\prod_{g\in G}g(s)\in S^G$. Hence $B=S^{-1}A=A\otimes_{A^G}(S^G)^{-1}A^G$. Since localisation is exact it follows that $B^G=(S^G)^{-1}A^G$.
\eproof

\subsection*{Calculation by resolutions}

We will use the following description of the graded ring $R_{A,G}=H^*(G,A)$. Let $P\to A$ be a resolution of $A$ by a complex $P$ of finite projective $AG$-modules, for example one can take a resolution $P'\to\ZZ$ by finite projective $\ZZ G$-modules and set $P=A\otimes_\ZZ P'$. Let 
\begin{equation}
\label{Eq:DefE}
E=\End_{AG}(P)
\end{equation}
as a differential graded algebra. Then 
\begin{equation}
\label{Eq:RAG-E}
R_{A,G}=H^*(E)
\end{equation}
as a graded ring. Moreover, let 
\begin{equation}
\label{Eq:DefN}
N=\Hom_{AG}(P,A)
\end{equation}
as a right dg $E$-module. The homomorphism $P\to A$ induces a quasi-isomorphism of right dg $E$-modules $E\to N$ and thus 
\begin{equation}
\label{Eq:RAG-N}
R_{A,G}=H^*(N)
\end{equation}
as a graded right $R_{A,G}$-module, in particular as a graded abelian group.

\bLe
\label{Le:Coh-uh-flat}
For a flat ring homomorphism $A^G\to B_0$ let $B=A\otimes_{A^G}B_0$ where $G$ acts on the first factor. Then $B^G=B_0$, and the homomorphism of graded rings $R_{A,G}\otimes_{A^G}B^G\to R_{B,G}$ is an isomorphism. In particular, the homomorphism $A\to B$ lies in the class $\Cohuh(G)$.
\eLe

\bproof
Flatness implies that $B^G=B_0$. For $N=\Hom_{AG}(P,A)$ as in \eqref{Eq:DefN} there is an isomorphism
\begin{equation}
\label{Eq:N-AG-BG}
N\otimes_{A^G}B^G\cong\Hom_{BG}(P\otimes_AB,B)
\end{equation}
because $P$ consists of finite projective $AG$-modules. Since $A^G\to B^G$ is flat, the cohomology of the left hand side of \eqref{Eq:N-AG-BG} is $R_{A,G}\otimes_{A^G}B^G$. Since $P\otimes_AB\to B$ is a resolution by finite projective $BG$-modules, the cohomology of the right hand side of \eqref{Eq:N-AG-BG} is $R_{B,G}$.
\eproof

\bPr
\label{Pr:Coh-uh-loc}
The class $\Cohuh(G)$ of Definition \ref{De:Coh-uh} contains all $G$-equivariant localisation homomorphisms.
\ePr

\bproof
By Lemma \ref{Le:AGBGloc} we have $B=A\otimes_{A^G}B_0$ where $B_0$ is a localisation of $A^G$, so the assertion follows from Lemma \ref{Le:Coh-uh-flat}.
\eproof

\bPr
\label{Pr:Coh-uh-nil}
Let $\bar A=A/I$ for a $G$-stable nilpotent ideal $I$. Then the natural ring homomorphism 
\[
\pi_{\even}:(R_{A,G})_{\even}\to (R_{\bar A,G})_{\even}
\] 
is a universal homeomorphism, and the projection $u:A\to\bar A$
lies in the class $\Cohuh(G)$ of Definition \ref{De:Coh-uh}.
\ePr

\bproof
For a prime $p$ let $p^r$ be the maximal $p$-power dividing $|G|$ and let $m_p=|G|/p^r$. Then $\Spec(\ZZ)$ is covered by the open sets $\Spec(\ZZ[1/m_p])$ for varying $p$, and we can replace $A$ by $A[1/m_p]$ for a fixed prime $p$, using Lemma~\ref{Le:Coh-uh-flat}. Then $p^rH^i(G,M)=0$ for any $AG$-module $M$ and $i>0$.

To prove that $\pi_{\even}$ is a universal homeomorphism, by an induction using the sequence of $G$-equivariant ring homomorphisms $A\to A/I^2\to A/I$ we can assume that $I$ has square zero. Then the exact sequence of $AG$-modules 
\[
0\to I\to A\to\bar A\to 0
\]
induces a long exact sequence in group cohomology
\begin{equation}
\label{Eq:HGIAA}
H^*(G,I)\xrightarrow j R_{A,G}\xrightarrow\pi R_{\bar A,G}\xrightarrow\delta H^*(G,I).
\end{equation}
Here $\pi$ is a homomorphism of graded rings that restricts to the homomorphism $\pi_{\even}$ of the proposition, $H^*(G,I)$ is a graded left and right $R_{\bar A,G}$-module since $I$ is an $\bar A$-module, and $j$ is $R_{A,G}$-linear. Moreover $\delta$ ist a graded derivation by Lemma \ref{Le:delta-deriv} below. If $a\in R_{\bar A,G}$ is homogeneous of even degree, we obtain
\begin{equation}
\label{Eq:delta-apr}
\delta(a^{p^r})=\sum_{i=1}^{p^r}a^{i-1}\delta(a)a^{p^r-i}=p^r\delta(a)a^{p^r-1}=0,
\end{equation}
hence $a^{p^{r}}$ lies in the image of $\pi$, moreover $p^ra$ lies in the image of $\pi$ since $p^r\delta(a)=0$. The image $j$ has square zero since this holds for $I$. It follows that $\pi_{\even}$ is a universal homeomorphism by \cite[\href{https://stacks.math.columbia.edu/tag/0BRA}{Lemma 0BRA}]{stacks-project}. 

To prove that $u$ lies in $\Cohuh(G)$ we can drop the assumption that $I$ has square zero. The degree zero component of $\pi_{\even}$ factors into the homomorphisms $A^G\to A^G/I^G\to \bar A^G$, which are both  universal homeomorphisms by Proposition \ref{Pr:AGIGBG} and since $I^G$ is nilpotent. Now $\pi_{\even}$ is the composition
\[
(R_{A,G})_{\even}\to (R_{A,G})_{\even}\otimes_{A^G}\bar A^G\xrightarrow{u'''} (R_{\bar A,G})_{\even}
\]
where the first arrow is a universal homeomorphism since this holds for $A^G\to\bar A^G$, and $u'''$ is the homomorphism associated to $u:A\to\bar A$ as in \eqref{Eq:RAGevRBGev}. Hence $u'''$ is a universal homeomorphism since this holds for $\pi_{\even}$.
\eproof

\bLe
\label{Le:delta-deriv}
The homomorphism $\delta$ of \eqref{Eq:HGIAA} is a graded derivation, i.e.\ $\delta(ab)=\delta(a)b+(-1)^{|a|}a\delta(b)$ for homogeneous elements $a,b\in R_{\bar A,G}$. 
\eLe

This is a variant of the well-known fact that the Bockstein homomorphism is a graded derivation, see for example \cite[Example 3E.1]{Hatcher:AT}. We include a proof for completeness.

\bproof
We use that $R_{A,G}=H^*(E)$ with $E=\End_{AG}(P)$ as in \eqref{Eq:DefE}. Let $\bar P=P/IP$. Then $\bar P\to\bar A$ is a resolution of $\bar A$ by finite projective $\bar AG$-modules, and $IP\to I$ is a resolution of $I$ by $\bar AG$-modules. Let $\bar E=\End_{\bar AG}(\bar P)$ and $J=\Hom_{AG}(P,IP)=\Hom_{\bar AG}(\bar P,IP)$. The obvious exact sequence $0\to IP\to P\to\bar P\to 0$ induces an exact sequence
\begin{equation}
\label{Eq:JEE}
0\to J\to E\xrightarrow{\tilde\pi}\bar E\to 0
\end{equation}
where $\tilde\pi$ is a homomorphism of dg algebras, so $J$ is a two-sided dg ideal of $E$, and $J$ becomes a dg $\bar E$-bimodule since $J$ has square zero. The cohomology sequence of \eqref{Eq:JEE} can be identified with \eqref{Eq:HGIAA}. For $a\in R_{\bar A,G}$ let $\tilde a\in E$ be an inverse image under $\tilde\pi$ of a representative of $a$ in $\bar E$. Then $\delta(a)=[d(\tilde a)]$ where $[\;]$ is the cohomology class of a cycle in $J$, and the lemma follows from the relation $d(\tilde a\tilde b)=d(\tilde a)\tilde b+(-1)^{|a|}\tilde ad(\tilde b)$.
\eproof

We recall that the pair $(G,A)$ is called noetherian if the ring $A^G$ is noetherian and if $A$ is finite over $A^G$; see Definition \ref{De:AG-noeth}.

\bPr
\label{Pr:Coh-uh-t}
Let $\bar A=A/tA$ where $t\in A^G$ is an $A$-regular element and assume that $(G,A)$ is noetherian. Then the natural ring homomorphism
\[
\bar\pi_{\even}:(R_{A,G}/tR_{A,G})_{\even}\to(R_{\bar A,G})_{\even}
\]
is a universal homeomorphism, and the projection $u:A\to\bar A$ lies in the class $\Cohuh(G)$ of Definition \ref{De:Coh-uh}.
\ePr

\bproof
We begin with two initial remarks. 

First, we can replace $t$ by a positive power $t^m$ using Proposition \ref{Pr:Coh-uh-nil} for to $A/t^m\to A/t$; note that the projection $(R_{A,G}/t^m)_{\even}\to(R_{A,G}/t)_{\even}$ is a universal homeomorphism, and $\Cohuh(G)$ is stable under composition by Proposition \ref{Pr:Coh-uh-comp}. The value of $m$ will be determined later.

Second, as in the proof of Proposition \ref{Pr:Coh-uh-nil}, after replacing $A$ by a localisation we can assume that for a fixed prime $p$ we have $p^rH^i(G,M)=0$ for any $AG$-module $M$ and $i>0$. 

The exact sequence of $AG$-modules 
\[
0\to A\xrightarrow tA\to \bar A\to 0
\]
induces a long exact sequence in group cohomology,
\begin{equation}
\label{Eq:RAGtRAG}
R_{A,G}\xrightarrow tR_{A,G}\xrightarrow\pi R_{\bar A,G}\xrightarrow\delta R_{A,G}\xrightarrow tR_{A,G}
\end{equation}
and thus an injective homomorphism of graded rings 
\begin{equation*}
\bar\pi:R_{A,G}/tR_{A,G}\to R_{\bar A,G}
\end{equation*}
that restricts to the homomorphism $\bar\pi_{\even}$ of the proposition. By Lemma \ref{Le:delta-apr} below, after replacing $t$ by $t^m$ for some positive $m$, for each homogeneous element $a\in R_{\bar A,G}$ of even degree, $a^{p^r}$ lies in the image of $\pi$. Moreover, $p^ra$ lies in the image of $\pi$ since $p^r\delta(a)=0$. It follows that $\bar\pi_{\even}$ is a universal homeomorphism by \cite[\href{https://stacks.math.columbia.edu/tag/0BRA}{Lemma 0BRA}]{stacks-project}. 

The degree zero component of $\bar\pi_{\even}$ is given by $A^G/tA^G\to \bar A^G$, which is a universal homeomorphism by Proposition \ref{Pr:AGIGBG}; note that $tA^G=(tA)^G$ since $t$ is $A$-regular. Now $\bar\pi_{\even}$ is the composition
\[
(R_{A,G})_{\even}\otimes_{A^G}A^G/tA^G\to (R_{A,G})_{\even}\otimes_{A^G}\bar A^G\xrightarrow{u'''} (R_{\bar A,G})_{\even}
\]
where the first arrow is a universal homeomorphism since this holds for $A^G/tA^G\to\bar A^G$, and $u'''$ is the homomorphism associated to $u:A\to\bar A$ as in \eqref{Eq:RAGevRBGev}. Hence $u'''$ is a universal homeomorphism since this holds for $\bar\pi_{\even}$.
\eproof

\bLe
\label{Le:delta-apr}
In the situation of Proposition \ref{Pr:Coh-uh-t}, after replacing $t$ by $t^m$ for some fixed positive integer $m$, for every homogeneous element $a\in R_{\bar A,G}$ of even degree we have $\delta(a^{p^{r}})=0$ in $R_{A,G}$.
\eLe

\bproof
This is similar to the calculation \eqref{Eq:delta-apr}, but with some complications since we do not have a direct analogue of Lemma \ref{Le:delta-deriv}. Since the ring $R_{A,G}$ is noetherian by Proposition \ref{Pr:RAG-noeth}, the ideals $J_i=\Ker(t^i:R_{A,G}\to R_{A,G})$ stabilize for large $i$. After replacing $t$ by $t^m$ for some $m\ge 1$ we can assume that $J_1=J_2$. We use again that $R_{A,G}=H^*(E)$ with $E=\End_{AG}(P)$ as in \eqref{Eq:RAG-E}. Let $\bar P=P/tP$. Then $\bar P\to\bar A$ is a resolution of $\bar A$ by finite projective $\bar AG$-modules. Let $\bar E=\End_{\bar AG}(\bar P)$. There is an exact sequence
\[
0\to E\xrightarrow t E\xrightarrow{\tilde\pi}\bar E\to 0,
\]
whose cohomology sequence can be identified with \eqref{Eq:RAGtRAG}. If $\tilde a\in E$ is an inverse image of a representative of $a$ in $\bar E$, then $\delta(a)=[d(\tilde a)/t]$ where $[\;]$ denotes the class of a cycle in $E$. Since $R_{\bar A,G}$ is a graded-commutative ring, for homogeneous elements $x,y,z\in E$ with $d(x),d(y),d(z)\in tE$ we have 
\begin{equation}
\label{Eq:xyz}
xyz-(-1)^{|x|\cdot|y|} yxz\in d(E)+tE.
\end{equation}
Modulo $d(E)+tE$ we obtain
\[
\frac{d(\tilde a^{p^r})}t=\sum_{i=1}^{p^r}\tilde a^{i-1}\frac{d(\tilde a)}t\tilde a^{p^r-i}\equiv p^r\frac{d(\tilde a)}t\tilde a^{p^r-1}\equiv 0,
\]
so $\delta(a^{p^r})=[d(x)+ty]=[ty]$ for certain $x,y\in E$. Necessarily, $d(x)+ty$ is a cycle, so $0=d(ty)=td(y)$. Hence $d(y)=0$ since $t$ is $E$-regular, and thus $\delta(a^{p^r})=tc$ with $c=[y]\in R_{A,G}$. We have $t^2c=t\delta(a^{p^r})=0$ since $t\circ\delta=0$ in \eqref{Eq:RAGtRAG}. Since $J_1=J_2$ it follows that $\delta(a^{p^r})=tc=0$.
\eproof

The following combination of Propositions \ref{Pr:Coh-uh-nil} and \ref{Pr:Coh-uh-t} will be used for a noetherian induction.

\bPr
\label{Pr:Coh-uh-AAq-prepare}
Assume that the pair $(G,A)$ is noetherian and the ring $A^G$ is local with maximal ideal $\Fq$. If the homomorphism $\psi_\Fq:A\to A(\Fq)$ of \eqref{Eq:psi_q} is not an isomorphism, then there is a non-zero $G$-invariant ideal $I$ of $A$ with $I\subseteq\Ker(\psi_\Fq)$ such that the projection $A\to A/I$ lies in the class $\Cohuh(G)$ of Definition \ref{De:Coh-uh}. 
\ePr

\bproof
If $A$ is not reduced, let $I\subseteq A$ be a non-zero $G$-stable nilpotent ideal, for example the nil-radical of $A$. Then $I\subseteq\Ker(\psi_\Fq)$ since $A(\Fq)$ is reduced, and the projection $A\to A/I$ lies in $\Cohuh(G)$ by Proposition \ref{Pr:Coh-uh-nil}. So we can assume that $A$ is reduced. 

By Lemma \ref{Le:AG-local}, the homomorphism $\psi_q$ is surjective and induces a homeomorphism $\Spec(A(\Fq))\cong\Max(A)$. If some maximal ideal of $A$ is also a minimal prime ideal, this holds for every maximal ideal of $A$ because $\Max(A)$ is a single $G$-orbit in $\Spec(A)$, hence $\Spec(A)=\Max(A)$, and $\psi_\Fq$ induces a bijective map $\Spec(A(\Fq))\to\Spec(A)$. Since $A$ is reduced and $\psi_\Fq$ surjective, this implies that $\psi_\Fq$ is bijective, which was excluded. Hence every minimal prime ideal $\Fp$ of $A$ is non-maximal and therefore satisfies $\Fp\cap A^G\ne\Fq$. 

By prime avoidance in $A^G$ we find an element $t\in A^G$ with $t\in\Fq$ such that $t$ is not contained in any minimal prime ideal of $A$; the minimal prime ideals of $A$ form a finite set since $A$ is noetherian. Since $A$ is reduced, $t$ is $A$-regular, so $A\to A/tA$ lies in $\Cohuh(G)$ by Proposition \ref{Pr:Coh-uh-t}.
\eproof

\bPr
\label{Pr:Coh-uh-AAq}
If the pair $(G,A)$ is noetherian, for each $\Fq\in\Spec(A^G)$ the homomorphism $\psi_\Fq:A\to A(\Fq)$ of \eqref{Eq:psi_q} lies in the class $\Cohuh(G)$ of Definition \ref{De:Coh-uh}.
\ePr

\bproof
Using the factorisation $A\to A_\Fq\to A(\Fq)$ of $\psi_\Fq$ of Lemma \ref{Le:psiq-factors}, Propositions \ref{Pr:Coh-uh-comp} and \ref{Pr:Coh-uh-loc} allow to replace $A$ by $A_\Fq$. Then $A^G$ is local with maximal ideal $\Fq$. The case $A=A(\Fq)$ is clear, so let $\psi_\Fq$ be not an isomorphism. Proposition \ref{Pr:Coh-uh-AAq-prepare} gives a factorisation of $\psi_\Fq$ into $G$-equivariant homomorphisms 
\begin{equation}
\label{Eq:Coh-uh-psiq-factors}
A\xrightarrow{\;\pi\;} A'=A/I\xrightarrow{\;\varphi\;} A(\Fq)
\end{equation}
where $I$ is non-zero and $\pi$ lies in $\Cohuh(G)$. 

Here $A'^G$ is local with maximal ideal $\Fq'$ lying over $\Fq$ because the natural map $\Spec(A'^G)\to\Spec(A^G/I^G)$ is a homeomorphism by Proposition \ref{Pr:AGIGBG}. The resulting homomorphism $A(\Fq)\to A'(\Fq')$ is an isomorphism by Lemma \ref{Le:AqBq-surj}, so $\varphi$ can be identified with the homomorphism \eqref{Eq:psi_q} for $A'$ and $\Fq'$ in place of $A$ and $\Fq$. One verifies that the pair $(G,A')$ is noetherian using the chain $A^G\to A'^G\subseteq A'$ where $A'$ is finite over $A^G$ and $A^G$ is noetherian.

Hence the hypotheses of Proposition \ref{Pr:Coh-uh-AAq-prepare} are satisfied by $A'$ and $\Fq'$ in place of $A$ and $\Fq$, and we can apply Proposition \ref{Pr:Coh-uh-AAq-prepare} repeatedly as long as the new rings $A'$ differ from $A(\Fq)$. Since $A$ is noetherian, the process necessarily stops and we arrive at a finite chain of $G$-equivariant ring homomorphisms $A\to A'\to A''\to\ldots\to A^{(n)}=A(\Fq)$ where all arrows lie in the class $\Cohuh(G)$. Then $\psi_\Fq$ lies in $\Cohuh(G)$ by Proposition \ref{Pr:Coh-uh-comp}.
\eproof

\bTh
\label{Th:univ-homeo}
The class $\Cohuh(G)$ of Definition \ref{De:Coh-uh} contains all $G$-equi\-variant ring homomorphisms $u:A\to B$ such that $B$ is a localisation of a quotient of $A$. 
\eTh

\bproof
The homomorphism $u:A\to B$ factors into $G$-equivariant homomorphisms $A\xrightarrow w A'\xrightarrow v B$ where $w$ is surjective and $v$ is an injective localisation. By Propositions \ref{Pr:Coh-uh-comp} and \ref{Pr:Coh-uh-loc} we can assume that $u$ is surjective. 

Let us first assume that the pair $(G,A)$ is noetherian in the sense of Definition \ref{De:AG-noeth}. Then $R_{B,G}$ is a finite module over the noetherian ring $R_{A,G}$ by Proposition \ref{Pr:RAG-noeth}, moreover $R_{A,G}$ is a finite module over the noetherian ring $(R_{A,G})_{\even}$ by Lemma \ref{Le:graded-noetherian}, and hence $(R_{B,G})_{\even}$ is finite over $(R_{A,G})_{\even}$. So the ring homomorphism $u'''$ of \eqref{Eq:RAGevRBGev} is finite, thus universally closed, and it suffices to show that $u'''$ is universally bijective. This holds iff for each $\Fq\in\Spec(B^G)$ the base change of $u'''$ under the natural ring homomorphism $B^G\to k(\Fq)$ is universally bijective, which is verified as follows.

Let $\Fq'\in\Spec(A^G)$ be the image of $\Fq$ and let $C=B(\Fq)$. Since  $A(\Fq')\cong B(\Fq)$ by Lemma \ref{Le:AqBq-surj}, in the sequence of $G$-equivariant homomorphisms $A\xrightarrow u B\to C$ the homomorphisms  $A\to C$ and $B\to C$ lie in $\Cohuh(G)$ by Proposition \ref{Pr:Coh-uh-AAq}. By two-out-of-three for universal homeomorphisms it follows that the base change of $u'''$ under the ring homomorphism $B^G\to C^G$ is a universal homeomorphism; see the proof of Proposition \ref{Pr:Coh-uh-comp}. By Lemma \ref{Le:AG-AqG-factors}, $B^G\to C^G$ factors as $B^G\to k(\Fq)\to C^G$ where $k(\Fq)\to C^G$ is a purely inseparable field extension and hence a universal homeomorphism. It follows that the base change of $u'''$ under under $B^G\to k(\Fq)$ is a universal homeomorphism as well. This finishes the proof if the pair $(G,A)$ is noetherian.

The general case follows by a limit argument. We have $A=\varinjlim_iA_i$ as a filtered direct limit where $A_i$ runs through all finitely generated $G$-invariant subrings of $A$, and $B=\varinjlim_iB_i$ where $B_i$ is the image of $A_i$ in $B$. Let $u_i:A_i\to B_i$ be the restriction of $u$. There is a commutative diagram
\[
\vcenter{
\xymatrix@C+1em@M+0.2em{
\varinjlim_i((R_{A_i,G})_{\even}\otimes_{A_i^G}B_i^G) \ar[r]^-{\varinjlim u_i'''} \ar[d] &
\varinjlim_i (R_{B_i,G})_{\even} \ar[d] \\
(R_{A,G})_{\even}\otimes_{A^G}B^G \ar[r]^-{u'''} &
(R_{B,G})_{\even}.
}}
\]
Each pair $(G,A_i)$ is noetherian by Corollary \ref{Co:AG-ft-noeth}, so $u_i$ lies in $\Cohuh(G)$ by the first part of the proof, i.e.\ $u_i'''$ is a universal homeomorphism. It follows that $\varinjlim u_i'''$ is a universal homeomorphism since $\Spec$ transforms a filtered colimit of rings into a limit of topological spaces; see \cite[Prop.~10.53]{GW:AG1}. The vertical arrows of the diagram are isomorphisms since the cohomology of a finite group preserves filtered colimits of the coefficients. Hence $u'''$ is a universal homeomorphism as desired.
\eproof

\bCo
\label{Co:Spech-Spec-cart}
Let $u:A\to B$ be a $G$-equivariant ring homomorphism where $B$ is a localisation of a quotient of $A$. Then there is a cartesian diagram of topological spaces with immersions as vertical arrows:
\begin{equation}
\label{Dia:Spech-Spc-cart}
\vcenter{
\xymatrix@M+0.2em{
\Spec^h(R_{B,G})
\ar[d]_{f^R} \ar[r]^-{\pi^R} &
\Spec(B^G) 
\ar[d]^{f^0} \\
\Spec^h(R_{A,G})
\ar[r]^-{\pi^R} &
\Spec(A^G)
}}
\end{equation}
\eCo

\bproof
The diagram is the right hand square of \eqref{Dia:TBH-TAG-extended} in the case $H=G$ and $f=u$. The homomorphism $u$ factors into a $G$-equivariant surjection and a $G$-equivariant localisation, and it suffices to treat these cases separately.

If $B$ is a $G$-equivariant localisation of $A$, then $B=S^{-1}A$ for a multiplicative set $S\subseteq A^G$ by Lemma \ref{Le:AGBGloc}, so we have $R_{B,G}=S^{-1}R_{A,G}$ by Lemma \ref{Le:Coh-uh-flat}, and all assertions follow. If $B=A/I$ for a $G$-invariant ideal $I$, then $A^G/I^G\to B^G$ is a universal homeomorphism by Proposition \ref{Pr:AGIGBG}. Hence in the chain of homomorphisms
\[
R_{A,G}\otimes_{A^G}A^G/I^G
\to
R_{A,G}\otimes_{A^G}B^G
\to
R_{B,G}
\]
the first arrow is a universal homeomorphism, while the second arrow is a universal homeomorphism by Theorem \ref{Th:univ-homeo}. Therefore, using Lemma \ref{Le:fh-surj-homeo}, the upper line of \eqref{Dia:Spech-Spc-cart} can be replaced by
\[
\Spec^h(R_{A,G}/I^GR_{A,G})\to\Spec(A^G/I^G),
\]
and again all assertions follow.
\eproof

\bCo
\label{Co:psiRq-homeo}
For an arbitrary pair $(G,A)$ and\/ $\Fq\in\Spec(A^G)$ the map
\[
\psi_{\Fq,\res}^R:\Spec^h(R_{A(\Fq),G})\to\Spec^h(R_{A,G})_\Fq
\]
in the fiber diagram \eqref{Dia:fiber-diagram} is a homeomorphism.
\eCo

\bproof
This follows from Corollary \ref{Co:Spech-Spec-cart} with $B=A(\Fq)$ because the continuous map $\Spec(A(\Fq)^G)\to\Spec(A^G)$ can be identified with the inclusion $\{\Fq\}\to\Spec(A^G)$ by Lemma \ref{Le:AG-AqG-factors}.
\eproof

\bRe
If the pair $(G,A)$ is noetherian, to deduce Corollary \ref{Co:psiRq-homeo} it is sufficient to use Proposition \ref{Pr:Coh-uh-AAq} instead of Theorem \ref{Th:univ-homeo}.
\eRe

\section{Tensor nilpotence}
\label{Se:TensorNil}

Let $F:\KKK\to\LLL$ be a tensor triangulated functor between essentially small tensor triangulated categories. Following \cite{Balmer:Surjectivity} we say that $F$ detects tensor nilpotence of morphisms if every morphism $f:X\to Y$ in $\KKK$ with $F(f)=0$ satisfies $f^{\otimes n}=0$ for some $n\ge 1$. 

If $F$ detects tensor nilpotence of morphisms and $\KKK$ is rigid, then the map $\Spc(F):\Spc(\LLL)\to\Spc(\KKK)$ is surjective by \cite[Thm.~1.3]{Balmer:Surjectivity}.

\bPr
\label{Pr:AA'-tensor-nilpotent}
Let $G$ be a finite group and let $f:A\to A'$ be a $G$-equivariant homomorphism of commutative rings. The resulting functor 
\[
f_*:\TTT_{A,G}\to\TTT_{A',G}
\] 
detects tensor nilpotence of morphisms in the following cases:
\begin{enumerate}
\item
\label{Eq:AA'-tensor-nilpotent-N}
$A'=A/N$ for a $G$-invariant nilpotent ideal $N$.
\item
\label{Eq:AA'-tensor-nilpotent-b}
$A'=A_b\times A/b$ for an $A$-regular element $b\in A^G$.
\end{enumerate}
\ePr

\bproof
We use the equivalence $\TTT_{A,G}\cong D^b(AG)_{A\pproj}$ given by Proposition \ref{Pr:TAG-DAGAperf} together with Lemma \ref{Le:DAG-Aproj-Aperf}. Let $f:X\to Y$ be a morphism in $D^b(AG)_{A\pproj}$ such that $f':X'\to Y'$ is zero in $D^b(A'G)$ where $X'=X\otimes_AA'$ etc. We have to show that $f^{\otimes n}=0$ for some $n\ge 1$.

We choose a quasi-isomorphism $u:P\to X$ where $P$ is a bounded above complex of finite projective $AG$-modules; see Lemma \ref{Le:DAGAperf-fin-res}. Then $f$ is represented by a homomorphism of complexes $g:P\to Y$ which is unique up to homotopy, namely $g=fu$ in $D(AG)$. The base change $g':P'\to Y'$ is homotopic to zero because $g'=f'u'$ in $D(A'G)$ where $f'=0$.

In the case \eqref{Eq:AA'-tensor-nilpotent-N} we write $g'=dh'+h'd$ for a homomorphism of graded $A'G$-modules $h':P'\to Y'[-1]$. Since $P$ consists of projective $AG$-modules, $h'$ lifts to a homomorphism of graded $AG$-modules $h:P\to Y[-1]$. We can replace $g$ by $g-(dh+hd)$ and thus assume that $g'$ is zero, which means that $g$ factors as $P\to NY\to Y$. Then $g^{\otimes r}$ factors as $P^{\otimes r}\to (NY)^{\otimes r}\to N^rY^{\otimes r}\to Y^{\otimes r}$. If $N^r=0$ it follows that $g^{\otimes r}=0$ as a homomorphism of complexes and thus $f^{\otimes r}=0$ in $D(AG)$. 
This refines \eqref{Eq:AA'-tensor-nilpotent-N} because the exponent $r$ is explicit (the nilpotence order of $N$).

In the case \eqref{Eq:AA'-tensor-nilpotent-b} let $X_b=X\otimes_AA_b$ etc.\  Since the components of $P$ are finite projective $AG$-modules we have
\begin{multline*}
\Hom_{D(A_bG)}(X_b,Y_b)=\Hom_{K(A_bG)}(P_b,Y_b)\\ {}=\Hom_{K(AG)}(P,Y)_b=\Hom_{D(AG)}(X,Y)_b.
\end{multline*}
The assumption $f'=0$ in $D(AG)$ implies that $f_b:X_b\to Y_b$ is zero in $D(A_bG)$, and hence $b^rf=0$ in $D(AG)$ for some $r\ge 1$. The obvious exact triangle $X\xrightarrow{b^r}X\xrightarrow\pi X/b^r\to^+$ in $D(AG)$ then shows that $f$ factors as
\[
X\xrightarrow\pi X/b^r\xrightarrow{\tilde f}Y.
\]
Now there is the following commutative diagram in $D(AG)$ where $(X/b^r)^{\otimes m}$ denotes the iterated tensor product over $A/b^r$, and similarly for $Y$, and where $f_r:X/b^r\to Y/b^r$ is the reduction of $f$.
Indeed, the upper square and triangle are evident, and the lower square follows from the corresponding square in $K(AG)$ with $P$ in place of $X$ and $g$ in place of $f$. 
\[
\vcenter{
\xymatrix@M+0.2em@C+1em{
X^{\otimes r}\otimes_A X
\ar[r]^{f^{\otimes r}\otimes{\id}} \ar[d]_{{\id}\otimes\pi} &
Y^{\otimes r}\otimes_A X
\ar[r]^{{\id}\otimes f} \ar[d]_{{\id}\otimes\pi} &
Y^{\otimes r}\otimes_A Y \\
X^{\otimes r}\otimes_A X/b^r
\ar[r]^{f^{\otimes r}\otimes\id}  \ar[d]_\simeq &
Y^{\otimes r}\otimes_A X/b^r
\ar[ur]_{{\id}\otimes\tilde f} \ar[d]_\simeq \\
(X/b^r)^{\otimes r}\otimes_{A/b^r}X/b^r 
\ar[r]^{f_r^{\otimes r}\otimes\id}  &
(Y/b^r)^{\otimes r}\otimes_{A/b^r}X/b^r \\
}}
\]
The assumption $f'=0$ implies that $f_1:X/b\to Y/b$ is zero in $D((A/b)G)$. By the refined version of \eqref{Eq:AA'-tensor-nilpotent-N} applied to the $G$-equivariant homomorphism $A/b^r\to A/b$ it follows that $f_r^{\otimes r}$ is zero in $D((A/b^r)G)$ and hence in $D(AG)$. It follows that $f^{\otimes (r+1)}=0$ in $D(AG)$.
\eproof

\section{Change of coefficients for $\protect\Spc$}
\label{Se:ChangeCoeffSpc}

As in \S \ref{Se:change-coefficient} we fix a finite group $G$.
For a $G$-equivariant ring homomorphism $f:A\to B$ we consider the outer square of the extended functoriality diagram \eqref{Dia:TBH-TAG-extended} for $G=H$:
\begin{equation}
\label{Dia:TBG-TAG-outer}
\vcenter{
\xymatrix@M+0.2em{
\Spc(\TTT_{B,G})
\ar[d]_{f^\TTT} \ar[r]^-{\bar\rho_{B,G}} &
\Spec(B^G) 
\ar[d]^{f^0} \\
\Spc(\TTT_{A,G})
\ar[r]^-{\bar\rho_{A,G}} &
\Spec(A^G)
}}
\end{equation}
and the corresponding fiber product of topological spaces
\[
(f^0)^*\Spc(\TTT_{A,G})=\Spc(\TTT_{A,G})\times_{\Spec(A^G)}\Spec(B^G).
\]
The diagram \eqref{Dia:TBG-TAG-outer} induces a continuous map
\begin{equation}
\label{Eq:f-TTT-res}
f^\TTT_{\res}:\Spc(\TTT_{B,G})\to (f^0)^*\Spc(\TTT_{A,G}),
\end{equation}
which we call the base change map associated to $f$. 

\bDe
\label{De:Spc-surj}
We denote by $\Spcsurj(G)$ the class of all $G$-equivariant homomorphisms of commutative rings $f:A\to B$ such that the base change map $f_{\res}^\TTT$ is surjective.
\eDe

One should compare this with Definition \ref{De:Coh-uh}. 

\bRe
\label{Re:rel-bc-map}
For the homomorphism $f=\psi_\Fq:A\to A(\Fq)$ of \eqref{Eq:psi_q} the base change map $f_{\res}^\TTT$ can be identified with the map $\psi_{\Fq,\res}^\TTT$ in the fiber diagram \eqref{Dia:fiber-diagram} because the natural map $\Spec(A(\Fq)^G)\to\Spec(A^G)$ can be identified with the inclusion $\{\Fq\}\to\Spec(A^G)$ by Lemma \ref{Le:AG-AqG-factors}. 
\eRe

\bPr
\label{Pr:Coh-surj-comp}
The class $\Spcsurj(G)$ is stable under composition.
\ePr

\bproof
If $A\xrightarrow{f}B\xrightarrow{g}C$ is a sequence of $G$-equivariant ring homomorphisms, then $(g\circ f)_{\res}^\TTT$ factors as
\[
\Spc(\TTT_{C,G})\xrightarrow{\;g_{\res}^\TTT\;}
(g^0)^*\Spc(\TTT_{B,G})
\xrightarrow{\;(g^0)^*f_{\res}^\TTT\;}
(g^0)^*(f^0)^*\Spc(\TTT_{A,G}).
\qedhere
\]
\eproof

\bPr
\label{Pr:Coh-surj-localis}
The class $\Spcsurj(G)$ contains all $G$-equivariant localisation homomorphisms.
\ePr

\bproof
Let $f:A\to B$ be a $G$-equivariant localisation, so $B=S^{-1}A$ for a multiplicative set $S\subseteq A^G$ by Lemma \ref{Le:AGBGloc}. By \cite[Cor.~3.10]{Balmer:Spectra3}, the Verdier localisation of $\TTT_{A,G}$ at the thick tensor ideal generated by $\cone(s)$ for $s\in S$ is a tensor triangulated category $S^{-1}\TTT_{A,G}$ with $\End(\One)=S^{-1}A^G=B^G$. Since the functor $f_*:\TTT_{A,G}\to\TTT_{B,G}$ maps $S$ to isomorphisms, it factors into tensor triangulated functors
\[
\TTT_{A,G}\xrightarrow{\;j\;} S^{-1}\TTT_{A,G}\xrightarrow{\;\varphi\;}\TTT_{B,G},
\]
and accordingly the transpose of the diagram \eqref{Dia:TBG-TAG-outer} factors as follows.
\begin{equation}
\label{Dia:TBG-TAG-outer-extend}
\vcenter{
\xymatrix@C+1em@M+0.2em{
\Spc(\TTT_{B,G})
\ar[r]^-{\Spc(\varphi)} \ar[d]_{\bar\rho_{B,G}} &
\Spec(S^{-1}\TTT_{A,G}) \ar[r]^-{\Spc(j)} \ar[d]^{\bar\rho} &
\Spc(\TTT_{A,G})
\ar[d]^{\bar\rho_{A,G}} \\
\Spec(B^G) \ar[r]^-{\sim} &
\Spec(S^{-1}A^G) \ar[r] &
\Spec(A^G)
}}
\end{equation}
Here the right hand square is cartesian by \cite[Thm.~5.4]{Balmer:Spectra3}.
The functor $\varphi$ is fully faithful because homomorphisms in $D^b(AG)_{A\pproj}$ are homomorphisms in $K^-(AG\pproj)$ by Lemma \ref{Le:DAGAperf-fin-res}, and these commute with localisation at $S$ by finiteness. Hence $\Spc(\varphi)$ is surjective by \cite[Cor.~1.8]{Balmer:Surjectivity}. Together if follows that $f_{\res}^\TTT$ is surjective.
\eproof

\bPr
\label{Pr:Coh-surj-nilp}
The class $\Spcsurj(G)$ contains all surjective $G$-equi\-variant homomorphisms $f:A\to B$ with nilpotent kernel.
\ePr

\bproof
Let $N$ be the kernel of $A\to B$. The map $f^0$ in \eqref{Dia:TBG-TAG-outer} factors into $\Spec(B^G)\to\Spec(A^G/N^G)\to\Spec(A^G)$ where both maps are homeomorphisms by Proposition \ref{Pr:AGIGBG} and because $N^G$ is a nilpotent ideal. So $f^0$ is a homeomorphism. The functor $f_*:\TTT_{A,G}\to\TTT_{B,G}$ detects tensor nilpotence by Proposition \ref{Pr:AA'-tensor-nilpotent} \eqref{Eq:AA'-tensor-nilpotent-N}, so the map $f^\TTT=\Spc(f_*)$ is surjective by \cite[Thm.~1.1]{Balmer:Surjectivity}. This map factors as
\begin{equation}
\label{Eq:Coh-surj-nilp-factors}
\Spc(\TTT_{B,G})
\xrightarrow{\;f_{\res}^\TTT\;}
(f^0)^*\Spc(\TTT_{A,G})
\xrightarrow{\;\pi\;}
\Spc(\TTT_{A,G})
\end{equation}
where $\pi$ is a base change of $f^0$ and thus bijective. Hence $f_{\res}^\TTT$ is surjective.
\eproof

\bPr
\label{Pr:Coh-surj-regular}
For each $A$-regular element $b\in A^G$, the projection homomorphism $g:A\to A/b$ lies in the class $\Spcsurj(G)$. 
\ePr

\bproof
Let $B=A_b\times A/b$. We consider the diagram \eqref{Dia:TBG-TAG-outer} for the natural $G$-equivariant homomorphism $f:A\to B$. Since $\TTT_{B,G}\cong\TTT_{A_b,G}\times\TTT_{A/b,G}$ as tensor triangulated categories and $B^G=(A_b)^G\times (A/b)^G$ as rings, there are compatible disjoint decompositions into open and closed subspaces
\begin{align}
\label{Eq:decomps-AB-Spc}
\Spc(\TTT_{B,G})&\cong\Spc(\TTT_{A_b,G})\sqcup\Spc(\TTT_{A/b,G}), 
\\
\label{Eq:decomps-AB-Spec}
\Spec(B^G)&\cong\Spec((A_b)^G)\sqcup\Spec((A/b)^G),
\end{align}
using Lemma \ref{Le:Spc-prod} for \eqref{Eq:decomps-AB-Spc}.
There is a similar disjoint decomposition into an open and a closed subspace
\begin{equation}
\label{Eq:decomp-Abb}
\Spec(A^G)\cong\Spec((A^G)_b)\sqcup\Spec(A^G/b).
\end{equation}
Here $(A^G)_b=(A_b)^G$ and $A^G/b=A^G/(bA)^G$ since $b$ is $A$-regular. Hence the decompositions \eqref{Eq:decomps-AB-Spec} and \eqref{Eq:decomp-Abb} together with Proposition \ref{Pr:AGIGBG} applied to $A\to A/b$ show that the map $f^0:\Spec(B^G)\to\Spec(A^G)$ is bijective.

Now the situation is similar to Proposition \ref{Pr:Coh-surj-nilp}. The functor $f_*:\TTT_{A,G}\to\TTT_{B,G}$ detects tensor nilpotence of morphisms by Proposition \ref{Pr:AA'-tensor-nilpotent} \eqref{Eq:AA'-tensor-nilpotent-b}, so the map $f^\TTT=\Spc(f_*)$ is surjective by \cite[Thm.~1.1]{Balmer:Surjectivity}. Again this map factors as \eqref{Eq:Coh-surj-nilp-factors} where $\pi$ is bijective since $f^0$ is bijective, so $f_{\res}^\TTT$ is surjective. The decompositions \eqref{Eq:decomps-AB-Spc} and \eqref{Eq:decomps-AB-Spec} yield that 
\[
g_{\res}^\TTT:\Spc(\TTT_{A/b,G})\to (g^0)^*\Spc(\TTT_{A,G})
\] 
is a retract of $f_{\res}^\TTT$, so $g_{\res}^\TTT$ is surjective as well.
\eproof

\bPr
\label{Pr:Coh-surj-AAq-prepare}
Assume that the ring $A$ is noetherian and the ring $A^G$ is local with maximal ideal $\Fq$. If the homomorphism $\psi_\Fq:A\to A(\Fq)$ of \eqref{Eq:psi_q} is not an isomorphism, then there is a non-zero $G$-invariant ideal $I$ of $A$ with $I\subseteq\Ker(\psi_\Fq)$ such that the projection $A\to A/I$ lies in the class $\Spcsurj(G)$ of Definition \ref{De:Spc-surj}. 
\ePr

\bproof
This is parallel to Proposition \ref{Pr:Coh-uh-AAq-prepare}, using Propositions \ref{Pr:Coh-surj-nilp} and \ref{Pr:Coh-surj-regular} instead of Proposition \ref{Pr:Coh-uh-nil} and \ref{Pr:Coh-uh-t}.
\eproof

\bPr
\label{Pr:Coh-surj-AAq}
If the ring $A$ is noetherian, then for each $\Fq\in\Spec(A^G)$ the homomorphism $\psi_\Fq:A\to A(\Fq)$ of \eqref{Eq:psi_q} lies in the class $\Spcsurj(G)$; in other words, the map $\psi_{\Fq,\res}^\TTT$ in the fiber diagram \eqref{Dia:fiber-diagram} is surjective.
\ePr

\bproof
This is parallel to Proposition \ref{Pr:Coh-uh-AAq}, using Proposition \ref{Pr:Coh-surj-comp}, \ref{Pr:Coh-surj-localis}, and \ref{Pr:Coh-surj-AAq-prepare} instead of Propositions \ref{Pr:Coh-uh-comp}, \ref{Pr:Coh-uh-loc}, and \ref{Pr:Coh-uh-AAq-prepare}. See Remark \ref{Re:rel-bc-map} for the assertion `in other words'.
\eproof

\section{The comparison map: conclusion}
\label{Se:main}

\bTh
\label{Th:rho-homeo}
For every pair $(G,A)$ the map 
\[
\rho_{A,G}:\Spc(\TTT_{A,G})\to\Spec^h(R_{A,G})
\]
of \eqref{Eq:rhoAG} is a homeomorphism.
\eTh

\bproof
We recall that $\TTT_{A,G}$ is equivalent to $D^b(AG)_{A\pperf}$ by Proposition \ref{Pr:TAG-DAGAperf}. 
Since $A$ is the filtered colimit of its finitely generated $G$-invariant subrings, by Lemma \ref{Le:DAGAperf-colim} and Lemma \ref{Le:Spc-colim} we can assume that $A$ is of finite type, in particular the pair $(G,A)$ is noetherian in the sense of Definition \ref{De:AG-noeth} by Corollary \ref{Co:AG-ft-noeth}. Then the rigid tensor category $\TTT_{A,G}$ is End-finite by Proposition \ref{Pr:RAG-noeth}, so by Corollary \ref{Co:rho-bij-homeo} the map $\rho_{A,G}$ is a homeomorphism iff it is bijective. By Lemma \ref{Le:rhobij-rhoqbij} this holds iff for each $\Fq\in\Spec(A^G)$ the map $(\rho_{A,G})_\Fq$ of \eqref{Eq:rhoAGq} is bijective. This map is the lower arrow in the fiber diagram \eqref{Dia:fiber-diagram}, in which the upper arrow $\rho_{A(\Fq),G}$ is a homeomorphism by Lemma \ref{Le:TAqG-TLH} and Theorem \ref{Th:field}, the right arrow $\psi_{\Fq,\res}^R$ is a homeomorphism by Corollary \ref{Co:psiRq-homeo}, and the left arrow $\psi_{\Fq,\res}^\TTT$ is surjective by Proposition \ref{Pr:Coh-surj-AAq}. It follows that $(\rho_{A,G})_\Fq$ is bijective as required.
\eproof

\bRe
In the first version of this article, Theorem \ref{Th:rho-homeo} was proved only when the ring $A$ is regular by a similar route, including also results on the functoriality of cohomological support. In view of Proposition \ref{Pr:rho-homeo-supp} these results are now obsolete.
\eRe

\section{A stable variant}
\label{Se:stable}

Let us record a variant of Theorem \ref{Th:rho-homeo} with $\Proj$ in place of $\Spec^h$, which is a rather formal consequence. This was first observed in \cite[\S 3.4]{Barthel:StratifyingIntegral} when $A$ is regular with trivial action of $G$. We begin with an elementary remark.

\bRe
The restriction of scalars $\res:AG\MMod\to A\MMod$ has a left adjoint $\ind:A\MMod\to AG\MMod$ defined by $\ind(M)=AG\otimes_AM$, using the right $A$-module structure of $AG$ for the tensor product. Let $P_0=AG$ as a left $AG$-module. For an $AG$-module $Q$ there is a natural isomorphism
\begin{equation}
\label{Eq:indresQ}
\ind\res Q\cong P_0\otimes_AQ
\end{equation}
where the tensor product is formed as in \eqref{Eq:PxQPxAQ}, i.e.\ $G$ acts diagonally. Indeed, the $A$-linear map $\res Q\to P_0\otimes_AQ$, $x\mapsto 1\otimes x$ gives by adjunction an $AG$-linear map $\ind\res Q\to P_0\otimes_AQ$, which is an isomorphism.
\eRe

\bLe
\label{Le:PerfAG}
The category $\Perf(AG)$ of perfect complexes of $AG$-modules is a thick tensor ideal in $D^b(AG)_{A\pperf}$. 
\eLe

\bproof
Clearly $\Perf(AG)$ is a thick subcategory of $D^b(AG)$. By Lemma \ref{Le:DAG-Aproj-Aperf} it suffices to show that for a finite $AG$-modules $P$, $Q$ where $P$ is projective and $Q$ is $A$-projective, $P\otimes_AQ$ is $AG$-projective. We can assume that $P=AG$. Then $P\otimes_AQ\cong\ind\res Q$ by \eqref{Eq:indresQ}, which is $AG$-projective.
\eproof

Using Proposition \ref{Pr:TAG-DAGAperf} we identify $\TTT_{A,G}$ and $D^b(AG)_{A\pperf}$. Then $\Perf(AG)$ is a tensor ideal in $\TTT_{A,G}$ by Lemma \ref{Le:PerfAG}, so the Verdier quotient
\[
\SSS_{A,G}=\TTT_{A,G}/\Perf(AG)
\] 
is a tensor triangulated category. 

\bCo
For every pair $(G,A)$ the homeomorphism $\rho_{A,G}$ of Theorem \ref{Th:rho-homeo} restricts to a homeomorphism
\[
\bar\rho_{A,G}:\Spc(\SSS_{A,G})\cong\Proj(R_{A,G}).
\]
\eCo

\bproof
Cf.\ \cite[Cor.~3.32]{Barthel:StratifyingIntegral}.
The ring homomorphism $q:R_{A,G}\to A^G$ defined by projection to degree zero gives a closed immersion of topological spaces $\Spec(A^G)\to\Spec^h(R_{A,G})$, whose complement is $\Proj(R_{A,G})$; moreover the natural functor $\TTT_{A,G}\to\SSS_{A,G}$ induces a homeomorphism
\begin{equation*}
\label{Eq:SpcSAG}
\Spc(\SSS_{A,G})\cong\{\PPP\in\Spc(\TTT_{A,G})\mid\Perf(AG)\subseteq\PPP\}
\end{equation*}
by \cite[Prop.~3.11]{Balmer:Spectrum}.
Hence it suffices to show that $\rho_{A,G}$ induces by restriction a bijective map $\bar\rho_{A,G}$ as indicated. But some $\PPP\in\Spc(\TTT_{A,G})$ satisfies $\Perf(AG)\subseteq\PPP$ iff $AG\in\PPP$ iff $\PPP\not\in\supp(AG)$ iff $\rho(\PPP)\not\in V(AG)$ by Proposition \ref{Pr:rho-homeo-supp}. Now $\End^*(AG)=AG$ in degree zero, which is an $R_{A,G}$-module via $q$, and it follows that $V(AG)=\Spec(A^G)$. Hence $\rho(\PPP)\not\in V(AG)$ iff $\rho(\PPP)\in\Proj(R_{A,G})$.
\eproof

\bRe
\label{Re:stable}
The Verdier quotient $\SSS_{A,G}$ can be viewed as the stable category associated to the Frobenius category $\lat(A,G)$ of $A$-projec\-tive finite $AG$-modules. More precisely, $\lat(A,G)$ is an exact subcategory of the abelian category $AG\MMod$, and with this exact structure $\lat(A,G)$ is a Frobenius category where the projective objects are the finite projective $AG$-modules; note that exact sequences in $\lat(A,G)$ are automatically $A$-split, and $\lat(A,G)$ has a duality involution defined by $M^\vee=\Hom_A(M,A)$ with $G$-action by conjugation. The stable category $\underline\lat(A,G)$ is a triangulated category by \cite[Thm.~2.6]{Happel:Triangulated}, and it is a tensor triangulated category since the projective objects of $\lat(A,G)$ form a tensor ideal. By an obvious extension of \cite[Thm.~2.1]{Rickard:Derived-Stable} there is a tensor triangulated equivalence
\begin{equation}
\label{Eq:ulat-SSS}
\underline\lat(A,G)\cong\SSS_{A,G}.
\end{equation}
In more detail, the composition $\lat(A,G)\to D^b(AG)_{A\pperf}\to\SSS_{A,G}$ induces the functor \eqref{Eq:ulat-SSS}, and an inverse functor is given by stabilized syzygies as follows.
The category $D^b(AG)_{A\pperf}$ is equivalent to the homotopy category $\KKK=K^-(AG\pproj)_{A\pperf}$ of upper bounded complexes of finite projective $AG$-modules which are $A$-perfect. Let $T$ denote the suspension functor of $\underline\lat(AG)$. Each $X\in\KKK$ has bounded cohomology. If $X$ has trivial cohomology in degree $\le n$, then the $AG$-module $Z^n(X)=\ker(d:X^n\to X^{n+1})$ is $A$-projective and the object $T^{-n}Z^n(X)$ of $\underline\lat(A,G)$ is independent of $n$. This construction gives a functor $\KKK\to\underline\lat(A,G)$, which induces an inverse of \eqref{Eq:ulat-SSS} as is easily verified. If $G$ acts trivially on $A$, the equivalence \eqref{Eq:ulat-SSS} is also proved in \cite[Prop.~3.26]{Barthel:StratifyingIntegral} using homotopy theoretic arguments.
\eRe

\bRe
If $G$ acts trivially on $A$, a different stable category $\stmod(AG)$ is considered in \cite{BIK:ModuleCategories}: The category of all $AG$-modules with the exact structure given by the $A$-split exact sequences is a Frobenius category, the associated stable category is denoted by $\StMod(AG)$, and $\stmod(AG)$ is the full subcategory of $\StMod(AG)$ whose objects are the finitely presented $AG$-modules. There is an obivous functor $\SSS_{A,G}\cong\underline\lat(A,G)\to\stmod(AG)$, but these categories behave quite differently with respect to tensor ideals. For example, for $A=\ZZ$ and a prime $p$ dividing the order of $G$, by \cite[\S 7]{BIK:ModuleCategories} there is an infinite ascending sequence $D_1\subsetneq D_2\subsetneq\ldots$ of radical tensor ideals in $\stmod(AG)$, where $D_n$ consisits of all finite $AG$-modules isomorphic to modules annihilated by $p^n$. This sequence is not visible in $\underline\lat(A,G)$ because the inverse image of $D_n$ in $\underline\lat(A,G)$ is zero for all $n$.
\eRe

\bRe
We also note the following geometric description of $\SSS_{A,G}$. The natural morphism $h:\Spec A\to[\Spec(A)/G]$ gives a direct image functor $h_*:\TTT_A\to\TTT_{A,G}$, which corresponds to the functor 
\[
\Perf(A)\to D^b(AG),\qquad
Q\mapsto AG\otimes_AQ.
\]
It follows that $\Perf(AG)$ coincides with $\left<h_*\TTT_A\right>$, the thick subcategory of $\TTT_{A,G}$ generated by the image of $h_*$, and hence $\SSS_{A,G}\cong\TTT_{A,G}/\left<h_*\TTT_A\right>$. 
\eRe

\end{document}